\documentclass[3p,times]{elsarticle}
\usepackage{epstopdf}
\usepackage{bbm}
\usepackage{bm}
\usepackage{array} 
\usepackage{color}
\usepackage{xcolor}

\usepackage{graphicx}
\usepackage{psfrag}

\usepackage{amssymb}
\usepackage{amsmath}
\usepackage{dsfont}
\usepackage{rotating}
\usepackage{pgf,tikz}
\usetikzlibrary{arrows}

\definecolor{ttzzqq}{rgb}{0.2,0.6,0}

\definecolor{qqttcc}{rgb}{0,0.2,0.8}

\definecolor{qqttzz}{rgb}{0,0.2,0.6}
\definecolor{ffqqqq}{rgb}{1,0,0}
\definecolor{qqwuqq}{rgb}{0,0.39,0}
\definecolor{zzttqq}{rgb}{0.6,0.2,0}
\definecolor{qqqqff}{rgb}{0,0,1}
\definecolor{ttttqq}{rgb}{0.2,0.2,0}

\definecolor{qqwwtt}{rgb}{0,0.4,0.2}
\definecolor{ubqqys}{rgb}{0.29,0,0.51}

\definecolor{wwttqq}{rgb}{0.4,0.2,0}
\definecolor{uuuuuu}{rgb}{0.27,0.27,0.27}
\definecolor{qqzzff}{rgb}{0,0.6,1}
\definecolor{xdxdff}{rgb}{0.49,0.49,1}
\definecolor{ccwwqq}{rgb}{0.8,0.4,0}

\definecolor{ttqqqq}{rgb}{0.2,0,0}

\definecolor{qqzzcc}{rgb}{0,0.6,0.8}

\newcommand{\flink}[2]{./#1#2-eps-converted-to.pdf}

\newcommand{\diff}[2]{\frac{\partial {#1} }{\partial {#2} } }
\newcommand{\eref}[1]{$(\ref{#1})$}
\newcommand{\E}{\mathbf{E}}
\newcommand{\R}{\mathbb{R}}
\newcommand{\xx}{ \mathbf{x}}
\newcommand{\dxt}{ d\mathbf{x}\, dt}
\newcommand{\dx}{ d\mathbf{x}}
\renewcommand{\u}{\mathbf{u}}
\newcommand{\q}{\mathbf{q}}
\newcommand{\uh}{\hat{\mathbf{u}}}
\newcommand{\Q}{\mathbf{Q}}
\newcommand{\B}{\mathbf{B}}
\renewcommand{\S}{\mathbf{S}}

\usepackage{amssymb}
\journal{Journal of Computational Physics}

\begin{document}

\begin{frontmatter}



\title{A simple diffuse interface approach on adaptive Cartesian grids for the linear elastic wave equations with complex topography}   

\author[1]{Maurizio Tavelli\fnref{label1}}
\author[1]{Michael Dumbser \corref{corr1} \fnref{label2}}
\author[2]{Dominic Etienne Charrier \fnref{label3}}
\author[3]{Leonhard Rannabauer \fnref{label4}}
\author[2]{\\Tobias Weinzierl \fnref{label5}}
\author[3]{Michael Bader \fnref{label6}}
\address[1]{Department of Civil, Environmental and Mechanical Engineering, University of Trento, Via Mesiano 77, I-38123 Trento, Italy}
\address[2]{Department of Computer Science, University of Durham, Lower Mountjoy, South Road, Durham DH1 3LE, United Kingdom}
\address[3]{Department of Informatics, Technical University Munich (TUM), Boltzmannstr. 3, D-85748 Garching, Germany}
\fntext[label1]{\tt m.tavelli@unitn.it (M.~Tavelli)}
\fntext[label2]{\tt michael.dumbser@unitn.it (M.~Dumbser)}
\fntext[label3]{\tt dominic.e.charrier@durham.ac.uk (D.E.~Charrier)}
\fntext[label4]{\tt leo.rannabauer@tum.de (L.~Rannabauer)}
\fntext[label5]{\tt tobias.weinzierl@durham.ac.uk (T.~Weinzierl)}
\fntext[label6]{\tt bader@in.tum.de (M.~Bader)}

\begin{abstract}
In most classical approaches of computational geophysics for seismic wave propagation problems, complex surface topography is either accounted for by boundary-fitted unstructured meshes, or, where possible, by mapping the complex computational domain from physical space to a topologically simple domain in a reference coordinate system. However, all these conventional approaches  
face problems if the geometry of the problem becomes sufficiently complex. They either need a mesh generator to create unstructured boundary-fitted grids, which can 
become quite difficult and may require a lot of manual user interactions in order to obtain a high quality mesh, or they need the explicit computation of an appropriate mapping function 
from physical to reference coordinates. For sufficiently complex geometries such mappings may either not exist or their Jacobian could become close to singular. 
Furthermore, in both conventional approaches low quality grids will always lead to very small time steps due to the Courant-Friedrichs-Lewy (CFL) condition for explicit schemes. 
In this paper, we propose a completely different strategy that follows the ideas of the successful family of high resolution shock-capturing schemes, where discontinuities can actually be resolved 
anywhere on the grid, without having to fit them exactly. We address the problem of geometrically complex free surface boundary conditions for seismic wave propagation problems with a 
novel diffuse interface method (DIM) on adaptive Cartesian meshes (AMR) that consists in the introduction of a characteristic function $ 0 \leq \alpha \leq 1 $ which identifies the location of the 
solid medium and the surrounding air (or vacuum) and thus implicitly defines the location of the free surface boundary. Physically, $\alpha$ represents the volume fraction of the solid medium present 
in a control volume. Our new approach \textit{completely avoids} the problem of mesh generation, since all that is needed for the definition of the complex surface topography is to set a scalar color  function to unity inside the regions covered by the solid and to zero outside. The governing equations are derived from ideas typically used in the mathematical description of compressible 
multiphase flows. An analysis of the eigenvalues of the PDE system shows that the complexity of the geometry has no influence on the admissible time step size due to the CFL condition. 
The model reduces to the classical linear elasticity equations inside the solid medium where the gradients of  $\alpha$ are zero, while in the diffuse interface zone at the free surface boundary the  governing PDE system becomes nonlinear.  
We can prove that the solution of the Riemann problem with arbitrary data and a jump in $\alpha$ from unity to zero yields a Godunov-state at the interface that satisfies the free-surface boundary condition exactly,  i.e.~the normal stress components vanish. In the general case of an interface that is not aligned with the grid and which is not infinitely thin, a systematic study on the  distribution of the volume fraction function inside the interface and the sensitivity with respect to the thickness of the diffuse interface layer has been carried out. In order 
to reduce numerical dissipation, we use high order discontinuous Galerkin (DG) finite element schemes on adaptive AMR grids together with a high resolution shock capturing subcell finite volume (FV)  limiter in the diffuse interface region. We furthermore employ a little dissipative HLLEM Riemann solver, which is able to resolve the steady contact discontinuity associated with the volume  
fraction function and the spatially variable material parameters exactly. 
	We provide a large set of computational results in two and three space dimensions involving complex geometries where the physical interface is not aligned with the grid. For all test cases we 
	provide a quantitative comparison with classical approaches based on boundary-fitted unstructured meshes.  
\end{abstract}


\end{frontmatter}


\section{Introduction}
The numerical solution of linear elastic wave propagation is still a challenging task, especially when complex three-dimensional geometries are involved. 
In the past, a large number of numerical schemes has been proposed for the simulation of seismic wave propagation. Madariaga \cite{Madariaga1976} and Virieux \cite{Virieux1984, Virieux1986} introduced finite difference schemes for the simulation of SH and P-SV wave propagation. These schemes were then extended to higher order, see \cite{Levander1988}, three space dimensions \cite{Mora1989,moczo2002} and to anisotropic material \cite{igel1995, tessmer1995}. For finite difference-like methods on unstructured meshes we refer to the work of  Magnier et. al.~\cite{magnier1994} and K\"aser \& Igel \cite{kaser2001a,kaser2001b}.
There are also several applications in the context of finite volume (FV) schemes \cite{Wang2002a,Wang2002b,Wang2004a,Wang2004b,Tadi2004,Dormt1995}, which, however, were all limited to second order of accuracy in space and time. The first arbitrary high order ADER finite volume scheme for seismic wave propagation was introduced in \cite{dumbserkaeser06c}.
For real applications it is crucial that a numerical scheme is able to properly capture complex signals over long  distances and times. 
 In contrast to classical low order schemes, high order methods in space and time are able to better reproduce the time evolution of the solution. 
 A quantitative accuracy analysis of high order numerical schemes for linear elasticity , based on the misfit criteria developed in \cite{Moczo2006,Moczo2009}, can be found in \cite{Hermann2008,Moczo2010a}. 
Spectral finite element methods \cite{Patera1984} were successfully applied to linear elastic wave propagation in a well-known series of paper of Komatitsch and collaborators \cite{Priolo1994,Komatitsch1998,Seriani1998,Komatitsch1999,Komatitsch2002}. 
For Chebyshev spectral method methods for wave propagation we refer to the work of Tessmer et. al. \cite{Tessmer1994, tessmer1995} and Igel \cite{Igel1999}. 
For alternative developments in the framework of stabilized continuous finite elements applied to elastic and acoustic wave propagation we refer to \cite{gugubrian_waves,song2015nitsche,scovazzi2017velocity}.
Apart from wave propagation in the medium, also the proper representation of complex surface topography is a challenging task.   
For this purpose, several high order numerical schemes on unstructured meshes were introduced in the past. A series of explicit high order discontinuous Galerkin (DG) schemes for elastic wave propagation on unstructured meshes was proposed in \cite{gij1, gij2,  gij5,GroteDG, Antonietti1, Antonietti2}, while the concept of space-time  discontinuous Galerkin schemes, originally introduced and analyzed in  \cite{spacetimedg1,spacetimedg2,KlaijVanDerVegt,Rhebergen2012,Rhebergen2013,Balazsova1,Balazsova2} for computational fluid dynamics (CFD), was later also extended to linear elasticity in \cite{Antonietti3,Antonietti4,SIStag2017}. 
The space-time DG method used in \cite{SIStag2017} is based on the novel concept of \textit{staggered} discontinuous Galerkin finite element schemes,  which was introduced for CFD problems in  \cite{DumbserCasulli,2DSIUSW,2STINS,3DSIINS,3DSIINS,Fambri2016,3DSICNS,AMRDGSI}. In any case, all previous methods require a boundary-fitted mesh that properly represent the geometry of the physical problem to be solved. 
The generation of this mesh is in general a highly non-trivial task and usually requires the use of external mesh generation tools. 
Moreover, the mesh generation process in highly complex geometry can lead to very small elements with bad aspect ratio, so-called {\itshape sliver elements} \cite{Bern1992,Joe1995,Fleischmann1999}. This well known problem can often be avoided, but not always, see e.g.  \cite{Cheng2000,Edelsbrunner2002}. 
For explicit time discretization, sliver elements can only be treated at the aid of {\itshape local time stepping} (LTS), see, for example, \cite{gij5,Taub2009,GroteLTS1,GroteLTS2}, but at the moment only very few  schemes used in production codes employed in computational seismology currently support time-accurate local  time stepping. Alternatively, implicit schemes  like  \cite{SIStag2017} require the introduction of a proper preconditioner in order to limit the number of iterations  needed to solve the associated linear algebraic system. 
%Our work is inspired by the ideas used for the simulation of compressible multiphase flows, see \cite{BaerNunziato1986,SaurelAbgrall}.

The \textit{key idea} of this paper is therefore to \textit{completely avoid} the mesh generation problem associated with classical approaches used in computational seismology. This is achieved by extending the linear elastic wave equations via a characteristic (color) function $\alpha$, which is nothing else than the volume fraction of the solid medium, and which determines if a point $\xx$ is located inside the solid material ($\alpha(\xx)=1$) or outside ($\alpha(\xx)=0$).  
In this way the scalar parameter $\alpha$ simply determines the physical boundary through a diffuse interface zone, instead of using a boundary-fitted unstructured mesh. With this new approach, even very complex  geometries can be easily represented on adaptive Cartesian meshes. Furthermore, the introduction of the new parameter $\alpha$ does not  change the eigenvalues of the PDE system and therefore does not influence the time step restriction imposed by the CFL condition. 
Our new method is inspired by the work concerning the modeling and simulation of compressible multiphase flows, see \cite{BaerNunziato1986} and \cite{SaurelAbgrall,SaurelAbgrall2,Abgral2003,AbgrallSaurel}. It can also be 
interpreted as a special case of the more general symmetric hyperbolic and thermodynamically compatible model of  nonlinear hyperelasticity of Godunov \& Romenski and collaborators  \cite{GodunovRomenski72,Rom1998,PeshRom2014,DumbserPeshkov2016,DumbserPeshkov2017,Romenskii2017}. 

A diffuse interface approach has already been successfully applied to nonlinear compressible fluid-structure interaction problems in a series of papers \cite{Ndanou2015,FavrieGavrilyuk2012,FavrieGavrilyuk2009}, but the employed numerical methods were low order accurate in space and time and therefore not suitable for seismic wave propagation problems. Other applications of diffuse interface methods for compressible multi-phase flows can be found in \cite{Gavrilyuk2008,Saurel2009,Kapila2001}, but, to the best of our  knowledge, this is the first time that a diffuse interface approach is derived and validated for the seismic wave propagation based on the equations of linear elasticity.  
Within the present paper, we use high order accurate ADER-DG schemes on Cartesian meshes with adaptive mesh refinement  (AMR). The numerical method has already successfully been applied to other hyperbolic PDE  systems \cite{Zanotti2015,DumbserPeshkov2016}. The use of adaptive mesh refinement allows to increase the resolution locally where needed, especially close to the free surface or at internal material boundaries. 
To avoid spurious oscillations and to enforce nonlinear stability, we use a simple but very robust  \textit{a posteriori} subcell finite volume limiter \cite{Dumbser2014}. Here, a second order total variation diminishing (TVD)  scheme is adopted in the limited zones on a finer sub-grid in order to maintain accuracy. The idea of using an \textit{a posteriori} approach to limit high order schemes was first proposed by Clain, Diot and Loub\`ere within the so-called Multi-dimensional Optimal Order Detection (MOOD) paradigm in the  context of finite volume schemes, see \cite{Clain2011,Diot2012} for more details. 
Finally, in our numerical scheme we make use of the HLLEM Riemann solver introduced in \cite{HLLEM,NCP_HLLEM}, which is able to resolve the steady contact discontinuities associated with the   spatially variable material parameters $\lambda$ and $\mu$ (the Lam\'e constants), the mass density 
$\rho$ and the volume fraction $\alpha$. The numerical results presented later in this paper show that the
proposed methodology seems to be a valid alternative to existing approaches in computational seismology that are based on boundary-fitted structured or unstructured meshes. 

The rest of the paper is organized as follows: in Section \ref{sec_MatMod} we introduce the governing PDE of the new diffuse interface approach for linear elasticity. We also show the compatibility of our model with the free surface boundary condition in the case where $\alpha$ jumps from $1$ to $0$. In Section  \ref{sec:numerical-scheme} we briefly summarize the high order ADER-DG schemes used in this paper. In  Section \ref{sec:numerical-results} we show numerical results for a large set of test problems in two and three space  dimensions, also including a realistic 3D scenario with complex geometry given by real DTM data. Finally, in Section \ref{sec:conclusions} we give some concluding remarks and an outlook on future work. 

\section{Mathematical model}
\label{sec_MatMod}
%For all the tests shown in this section we use the complete algorithm involving high 
%order ADER-DG schemes, AMR and the finite volume sub-cell limiter. The detection criterion is based on the observed quantity $\varphi(Q)=\alpha$, which identifies the location of the boundary. Furthermore, we will always use a Riemann solver that does not introduce numerical dissipation for the steady contact waves associated with the material parameters and $\alpha$, such as HLLEM \cite{HLLEM} or the Osher method \cite{OsherNC}.
\noindent The equations of linear elasticity \cite{BedfordDrumheller} can be written as
\begin{eqnarray}
	\diff{}{t}\sigma_{xx}-\diff{}{x}u-\lambda \diff{}{y}v-\lambda \diff{}{z}w = S_{xx}, \nonumber \\
	\diff{}{t}\sigma_{yy}-\lambda \diff{}{x}u-(\lambda+2\mu) \diff{}{y}v-\lambda \diff{}{z}w = S_{yy}, \nonumber \\
	\diff{}{t}\sigma_{zz}-\lambda \diff{}{x}u-\lambda \diff{}{y}v-(\lambda+2\mu)\diff{}{z}w  = S_{zz}, \nonumber \\ 
	\diff{}{t}\sigma_{xy}-\mu \left(\diff{}{x}v+\diff{}{y}u \right) = S_{xy}, \nonumber \\
	\diff{}{t}\sigma_{yz}-\mu \left(\diff{}{z}v+\diff{}{y}w \right) = S_{yz}, \nonumber \\
	\diff{}{t}\sigma_{xz}-\mu \left(\diff{}{z}u+\diff{}{x}w \right) = S_{xz}, \nonumber \\
	\diff{}{t}\left(\rho u\right)-\diff{}{x}\sigma_{xx}-\diff{}{y}\sigma_{xy}-\diff{}{z}\sigma_{xz}=\rho S_{u}, \nonumber \\
	\diff{}{t}\left(\rho v\right)-\diff{}{x}\sigma_{xy}-\diff{}{y}\sigma_{yy}-\diff{}{z}\sigma_{yz}=\rho S_{v}, \nonumber \\
	\diff{}{t}\left(\rho w\right)-\diff{}{x}\sigma_{xz}-\diff{}{y}\sigma_{yz}-\diff{}{z}\sigma_{zz}=\rho S_{w}. 
\label{eq:1}
\end{eqnarray}
In more compact form the above system reads
\begin{eqnarray}
\diff{\bm{\sigma}}{t}-\mathbf{E}(\lambda, \mu)\cdot \nabla {\bm{v}}=\bm{S}_\sigma, \label{eq:1.1} \\
\diff{\rho\bm{v}}{t}-\nabla \cdot \bm{\sigma}=\rho \bm{S}_v,	\label{eq:1.2} 
\end{eqnarray}
where $\bm{v}=(u,v,w)$ is the velocity field, $\rho$ is the material density,  $\bm{S}_\rho$ and $\bm{S}_\sigma$ are volume sources, $\bm{\sigma}$ is the symmetric stress tensor, and $\E(\lambda,\mu)$ is the stiffness tensor that connects the strain tensor $\epsilon_{kl}$ to the stress tensor $\bm{\sigma}$ according to the Hooke law $\bm{\sigma}=\mathbf{E}\epsilon$.
The stress tensor $\bm{\sigma}$ is given by
\begin{eqnarray}
	\bm{\sigma}=\left(
	\begin{array}{ccc}
		\sigma_{xx} & \sigma_{xy} & \sigma_{xz} \\
		\sigma_{yx} & \sigma_{yy} & \sigma_{yz} \\
		\sigma_{zx} & \sigma_{zy} & \sigma_{zz} \\
	\end{array} \right)
\label{eq:2}
\end{eqnarray}
with the symmetry $\sigma_{ij}=\sigma_{ji}$. The normal stress components are $\sigma_{xx}, \sigma_{yy}$ and $\sigma_{zz}$, while the shear stress is represented by $\sigma_{xy}, \sigma_{yz}$ and $\sigma_{xz}$. The stress tensor $\bm{\sigma}$ can thus be written in terms of its six independent components $(\sigma_{xx},\sigma_{yy},\sigma_{zz},\sigma_{xy},\sigma_{yz},\sigma_{xz})$.
In the following we propose a new model that follows the ideas used in the simulation of compressible multiphase flows \cite{BaerNunziato1986,SaurelAbgrall,SaurelAbgrall2,AbgrallSaurel}. In order to do derive the model, we start from the Baer-Nunziato system of compressible multi-phase flows, where for the solid phase (index $s$) the pressure term has been appropriately replaced by the stress tensor $\boldsymbol{\sigma}_s$:  

\begin{eqnarray}
%\left\{
%\renewcommand{\arraystretch}{2}
%\begin{array}{l}
\label{ec.BN}
%
% First equation
%
\frac{\partial}{\partial t}\left(\alpha_s \rho_s \right)+\nabla\cdot\left(\alpha_s \rho_s \textbf{v}_s\right) &=& 0, \nonumber
\\
%
% Second equation
%
\frac{\partial}{\partial t}\left(\alpha_s \rho_s \textbf{v}_s\right)
+\nabla\cdot\left(\alpha_s \rho_s \textbf{v}_s \otimes \textbf{v}_s + \alpha_s \boldsymbol{\sigma}_s \right) 
-\boldsymbol{\sigma}_I \nabla \alpha_s &=&  \alpha_s \rho_s \mathbf{S}_{v,s}, 
\nonumber \\
%
%
% Third equation
%
\frac{\partial}{\partial t}\left(\alpha_s \rho_s E_s\right) 
+ \nabla \cdot \left( \alpha_s \rho_s E_s \textbf{v}_s + \alpha_s \boldsymbol{\sigma}_s  \mathbf{v}_s \right) - \boldsymbol{\sigma}_I \nabla \alpha_s \cdot \mathbf{v}_I &=& \alpha_s \rho_s \mathbf{S}_{v,s} \cdot \mathbf{v}_s,
\nonumber \\
%
% Fourth equation
%
\frac{\partial}{\partial t}\left(\alpha_g \rho_g \right)+\nabla\cdot\left(\alpha_g \rho_g \textbf{v}_g\right) &=& 0, \nonumber
\nonumber \\ 
%
% Fifth equation
%
\frac{\partial}{\partial t}\left(\alpha_g \rho_g \textbf{v}_g\right)
+\nabla\cdot\left(\alpha_g \rho_g \textbf{v}_g \otimes \textbf{v}_g + \alpha_g \boldsymbol{\sigma}_g \right) 
-\boldsymbol{\sigma}_g \nabla \alpha_g &=&  \alpha_g \rho_g \mathbf{S}_{v,g}, 
\nonumber \\
%
% Sixth equation
%
\frac{\partial}{\partial t}\left(\alpha_g \rho_g E_g \right)
+\nabla\cdot\left( \alpha_g \rho_g E_g \textbf{v}_g + \alpha_g \boldsymbol{\sigma}_g  \mathbf{v}_g \right) - \boldsymbol{\sigma}_I \nabla \alpha_g \cdot \mathbf{v}_I &=& \alpha_g \rho_g \mathbf{S}_{v,g} \cdot \mathbf{v}_g,
\nonumber \\
%
% Seventh equation
%
\frac{\partial}{\partial t}\alpha_s + \textbf{v}_{I}\nabla\alpha_s &=&  0.
%\end{array}\right.
\end{eqnarray}
Here index $s$ refers to the solid phase and index $g$ refers to the gas phase surrounding the solid; $\rho_k$ is the mass density and $E_k$ is the specific total energy of phase $k$,  $\mathbf{v}_k$ is the phase velocity, $\mathbf{v}_I$ is the so-called interface velocity and 
$\boldsymbol{\sigma}_I$ is the stress tensor at the interface, which is a generalization of the interface pressure used in standard BN models. We now make the following simplifying assumptions: i) The interface between solid and gas is not moving, hence 
$\mathbf{v}_I=0$. ii) The mass density of the gas phase is much smaller than the one of the solid
phase ($\rho_g \ll \rho_s$), hence the time evolution of the gas phase is not relevant 
for our purposes. Therefore, all evolution equations related to the gas phase can be neglected in 
the following, similar to the approach used in \cite{DIM2D,DIM3D} in the context of non-hydrostatic free  surface flow simulations based on a diffuse interface approach. To ease notation, the remaining index $s$ for the solid phase can be dropped. 
iii) We assume the density $\rho_s$ of the solid phase to be constant in time, so the related mass conservation equation can be neglected. iv) Furthermore, the stress tensor of the solid can be directly calculated via Hooke's law \eqref{eq:1.1}, so it is not necessary to evolve the total energy conservation law for the solid. v) The nonlinear convective term $\alpha_s \rho_s \mathbf{v}_s \otimes \mathbf{v}_s$, which is quadratic in the solid velocity, can be neglected, since the solid velocity is assumed to be small in the linear elasticity limit. vi) Last but not least, the free surface boundary condition at the interface between solid and surrounding gas leads to $\boldsymbol{\sigma}_s \cdot \nabla \alpha_s =0$. As a result of these simplifying assumptions, the \textit{reduced} governing PDE system of the 
new diffuse interface approach for linear elasticity in complex geometry reads:   
\begin{eqnarray}
\diff{\bm{\sigma}}{t} - \E(\lambda, \mu)\cdot \nabla {\bm{v}} = \bm{S}_\sigma,  \\
\diff{\alpha \rho \bm{v}}{t} - \nabla \cdot (\alpha \bm{\sigma})  = \alpha \rho \bm{S}_v,  \\ 
\diff{\alpha}{t} = 0.  
\end{eqnarray}
Since $\partial_t \rho = 0$, the previous equations are then rewritten as 
\begin{eqnarray}
\diff{\bm{\sigma}}{t} - \E(\lambda, \mu) \cdot \frac{1}{\alpha} \nabla {(\alpha\bm{v})} + \frac{1}{\alpha}\E(\lambda, \mu) \cdot {\bm{v}} \otimes \nabla \alpha = \bm{S}_\sigma, \label{eq:3.00}\\
\diff{\alpha \bm{v}}{t}-\frac{\alpha}{\rho}\nabla \cdot \bm{\sigma} - \frac{1}{\rho}\sigma\nabla \alpha=\bm{S}_v, \\ 
\label{eq:3.0}
\diff{\alpha}{t} = 0.
\label{eq:3}
\end{eqnarray}
Furthermore the following equations for the material parameters are added to the system: 
\begin{equation}
	\diff{\lambda}{t}=0,   \qquad 	\diff{\mu}{t}=0,  \qquad \diff{\rho}{t}=0. \label{eq:3.1}
\end{equation}
The same diffuse interface model can also be obtained by combining the nonlinear hyperelasticity equations of Godunov and Romenski \cite{GodunovRomenski72,Godunov:2003a,Godunov:1995a} with the compressible multi-phase model of Romenski et al. \cite{Rom1998,Rom2010}, assuming linear material behavior and neglecting nonlinear convective terms.   
System \eref{eq:3.00}-\eref{eq:3.1} is then rewritten in the following form: 
\begin{eqnarray}
\diff{\mathbf{Q}}{t} 
+ \mathbf{B}_{1}(\Q) \diff{\mathbf{Q}}{x}  
+ \mathbf{B}_{2}(\Q)  \diff{\mathbf{Q}}{y} 
+ \mathbf{B}_{3}(\Q)  \diff{\mathbf{Q}}{z} =\bm{S}(\xx,t),  
\label{eq:4}
\end{eqnarray}
where the three matrices $\mathbf{B}_{1}$, $\mathbf{B}_{2}$ and $\mathbf{B}_{3}$ are specified in Eqs. \eref{eq:8}-\eref{eq:10}. The vector $\mathbf{Q}$ is given by
\begin{eqnarray}
	\mathbf{Q} = \left( \sigma_{xx},\sigma_{yy},\sigma_{zz},\sigma_{xy},\sigma_{yz},\sigma_{xz},\alpha u,\alpha v,\alpha w, \lambda, \mu, \rho, \alpha \right)^\top,
\label{eq:11}
\end{eqnarray}
while the matrices $\mathbf{B}_{1}$, $\mathbf{B}_{2}$ and $\mathbf{B}_{3}$ read 
\begin{eqnarray}
	\mathbf{B}_1 = \left( 
	\begin{array}{ccccccccccccc}
		0 & 0& 0&0&0&0&-\frac{1}{\alpha}(\lambda+2\mu) & 0&0&0&0&0&\frac{1}{\alpha}(\lambda+2\mu) u \\
		0 & 0& 0&0&0&0&-\frac{1}{\alpha}\lambda & 0&0&0&0&0&\frac{1}{\alpha}\lambda u \\
		0 & 0& 0&0&0&0&-\frac{1}{\alpha}\lambda & 0&0&0&0&0&\frac{1}{\alpha}\lambda u \\
		0 & 0& 0&0&0&0&0 & -\frac{1}{\alpha}\mu &0&0&0&0&\frac{1}{\alpha}\mu v \\
		0 & 0& 0&0&0&0&0 &0 &0 &0&0&0&0 \\
		0 & 0& 0&0&0&0&0 &0 &-\frac{1}{\alpha}\mu &0&0&0&\frac{1}{\alpha}\mu w \\
		-\frac{\alpha}{\rho} & 0& 0&0&0&0&0 &0 &0 &0&0&0&-\frac{1}{\rho} \sigma_{xx} \\
		0 & 0& 0&-\frac{\alpha}{\rho}&0&0&0 &0 &0 &0&0&0&-\frac{1}{\rho} \sigma_{xy} \\
		0 & 0& 0&0&0&-\frac{\alpha}{\rho}&0 &0 &0 &0&0&0&-\frac{1}{\rho} \sigma_{xz} \\
		0 & 0& 0&0&0&0&0 &0 &0 &0&0&0&0 \\
		0 & 0& 0&0&0&0&0 &0 &0 &0&0&0&0 \\
		0 & 0& 0&0&0&0&0 &0 &0 &0&0&0&0 \\
		0 & 0& 0&0&0&0&0 &0 &0 &0&0&0&0 \\
	\end{array}
	\right),
\label{eq:8}
\end{eqnarray}
\begin{eqnarray}
	\mathbf{B}_2 = \left( 
	\begin{array}{ccccccccccccc}
	  0 & 0& 0&0&0&0&0 &-\frac{1}{\alpha}\lambda& 0&0&0&0&\frac{1}{\alpha}\lambda v \\
		0 & 0& 0&0&0&0& 0&-\frac{1}{\alpha}(\lambda+2\mu) &0&0&0&0&\frac{1}{\alpha}(\lambda+2\mu) v  \\
		0 & 0& 0&0&0&0 &0 &-\frac{1}{\alpha}\lambda& 0&0&0&0&\frac{1}{\alpha}\lambda v  \\
		0 & 0& 0&0&0&0 & -\frac{1}{\alpha}\mu&0 &0&0&0&0&\frac{1}{\alpha}\mu u  \\
		0 & 0& 0&0&0&0&0&0 & -\frac{1}{\alpha}\mu &0&0&0&\frac{1}{\alpha}\mu w \\
		0 & 0& 0&0&0&0&0 &0 &0 &0&0&0&0 \\
		0 &0&0& -\frac{\alpha}{\rho}& 0&0&0 &0 &0 &0&0&0&-\frac{1}{\rho} \sigma_{xy} \\
		0 & -\frac{\alpha}{\rho}& 0&0&0&0&0 &0 &0 &0&0&0&-\frac{1}{\rho} \sigma_{yy} \\
		0 &0&0&0& -\frac{\alpha}{\rho}& 0&0 &0 &0 &0&0&0&-\frac{1}{\rho} \sigma_{yz} \\
		0 & 0& 0&0&0&0&0 &0 &0 &0&0&0&0 \\
		0 & 0& 0&0&0&0&0 &0 &0 &0&0&0&0 \\
		0 & 0& 0&0&0&0&0 &0 &0 &0&0&0&0 \\
		0 & 0& 0&0&0&0&0 &0 &0 &0&0&0&0 \\
	\end{array}
	\right),
\label{eq:9}
\end{eqnarray}
\begin{eqnarray}
	\mathbf{B}_3=\left( 
	\begin{array}{ccccccccccccc}
	  0 & 0& 0&0&0&0&0&0&-\frac{1}{\alpha}\lambda & 0&0&0&\frac{1}{\alpha}\lambda w  \\
		0 & 0& 0&0&0&0&0&0&-\frac{1}{\alpha}\lambda & 0&0&0&\frac{1}{\alpha}\lambda w  \\
		0 & 0& 0&0&0&0&0&0&-\frac{1}{\alpha}(\lambda+2\mu) & 0&0&0&\frac{1}{\alpha}(\lambda+2\mu) w  \\
		0 & 0& 0&0&0&0&0 &0 &0 &0&0&0&0 \\
		0 & 0& 0&0&0&0&0 & -\frac{1}{\alpha}\mu &0&0&0&0&\frac{1}{\alpha}\mu v  \\
		0 & 0& 0&0&0&0 & -\frac{1}{\alpha}\mu&0 &0&0&0&0&\frac{1}{\alpha}\mu u  \\
		0& 0&0&0& 0&-\frac{\alpha}{\rho} &0 &0 &0 &0&0&0&-\frac{1}{\rho} \sigma_{xz} \\
		0& 0&0&0&-\frac{\alpha}{\rho} & 0&0 &0 &0 &0&0&0&-\frac{1}{\rho} \sigma_{yz} \\
		0& 0&-\frac{\alpha}{\rho} &0&0& 0&0 &0 &0 &0&0&0&-\frac{1}{\rho} \sigma_{zz} \\
		0 & 0& 0&0&0&0&0 &0 &0 &0&0&0&0 \\
		0 & 0& 0&0&0&0&0 &0 &0 &0&0&0&0 \\
		0 & 0& 0&0&0&0&0 &0 &0 &0&0&0&0 \\
		0 & 0& 0&0&0&0&0 &0 &0 &0&0&0&0 \\
	\end{array}
	\right).
\label{eq:10}
\end{eqnarray}
The eigenvalues associated with the matrix $\mathbf{B}_{1}$ are 
\begin{equation}
 \lambda_{1} = -c_p, \quad \lambda_{2,3} = -c_s, \quad \lambda_{4,5,6,7,8,9,10} = 0, \quad \lambda_{11,12}=+c_s, \quad \lambda_{13} = +c_p, 
\label{eqn.eval} 
\end{equation}
where
\begin{equation}
	c_p=\sqrt{\frac{\lambda+2\mu}{\rho}} \qquad \textnormal{ and } \qquad c_s=\frac{\mu}{\rho}
\end{equation}
are the $p-$ and $s-$ wave velocities, respectively.
The matrix of right eigenvectors of the matrix $\mathbf{B}_1$ as defined in \eref{eq:8} is given by 
\begin{eqnarray}
	\mathbf{R} = \left( 
	\begin{array}{ccccccccccccc}
	  \rho c_p^2 & 0& 0&0&0&0&0 &0 &0 &-\sigma_{xx}&0&0&\rho c_p^2 \\
		\rho(c_p^2-2c_s^2) & 0& 0&1&0&0&0 &0 &0 &0&0&0&\rho(c_p^2-2c_s^2) \\
		\rho(c_p^2-2c_s^2) & 0& 0&0&1&0&0 &0 &0 &0&0&0&\rho(c_p^2-2c_s^2) \\
		0 & \rho c_s^2& 0&0&0&0&0 &0 &0 &-\sigma_{xy}&0&\rho c_s^2&0 \\
		0 & 0& 0&0&0&1&0 &0 &0 &0&0&0&0 \\
		0 & 0&  \rho c_s^2&0&0&0&0 &0 &0 &-\sigma_{xz}& \rho c_s^2&0&0 \\
		c_p & 0& 0&0&0&0&0 &0 &0 & \alpha u &0&0&-c_p \\
		0 & c_s& 0&0&0&0&0 &0 &0 & \alpha v &0&-c_s &0 \\
		0 & 0& c_s&0&0&0&0 &0 &0 & \alpha w &-c_s &0&0 \\
		0 & 0& 0&0&0&0&0 &0 &1 &0&0&0&0 \\
		0 & 0& 0&0&0&0&0 &1 &0 &0&0&0&0 \\
		0 & 0& 0&0&0&0&1 &0 &0 &0&0&0&0 \\
		0 & 0& 0&0&0&0&0 &0 &0 &\alpha&0&0&0 \\
	\end{array}
	\right).
\label{eq:RM}
\end{eqnarray}
The expressions for the eigenvalues and eigenvectors of $\B_2$ and $\B_3$ are very similar and can be obtained
from those of $\B_1$ since the PDE system is rotationally invariant. For this reason, we do not give their
explicit expressions here.  
We now want to show that the proposed model satisfies the free surface boundary condition $\boldsymbol{\sigma} \cdot \mathbf{n}=0$ exactly when considering a Riemann problem that includes a jump of $\alpha$ from $\alpha^L=1$ to  $\alpha^R=0$. 

For this, consider the left and right state of a Riemann problem in the $x$-direction given by 
\begin{eqnarray}
\mathbf{Q}_L & = & (\sigma_{xx}^L,\sigma_{yy}^L,\sigma_{zz}^L,\sigma_{xy}^L,\sigma_{yz}^L,\sigma_{xz}^L,u^L,v^L,w^L,\lambda, \mu, \rho,1 ), \label{eqn.ql} \\
\mathbf{Q}_R & = & (\sigma_{xx}^R,\sigma_{yy}^R,\sigma_{zz}^R,\sigma_{xy}^R,\sigma_{yz}^R,\sigma_{xz}^R,0,0,0,\lambda, \mu, \rho,0 ) \label{eqn.qr}.  
\end{eqnarray}
By using a simple straight line segment path   
\begin{equation}
   \psi(s) = \mathbf{Q}_L + s \left( \mathbf{Q}_R - \mathbf{Q}_L \right), 
\end{equation} 
we can define a generalized Roe-averaged matrix $\tilde{\mathbf{B}}_1$ in $x$ direction according to  \cite{Castro2006,Pares2006,OsherNC} as follows:   
\begin{equation}
   \tilde{\mathbf{B}}_1 = \int \limits_0^1 \mathbf{B}_1(\psi(s)) ds. 
\end{equation} 
The exact solution of the linearized Riemann problem based on the Roe-averaged matrix $\tilde{\mathbf{B}}_1 = \tilde{\mathbf{R}} \tilde{\boldsymbol{\Lambda}} \tilde{\mathbf{R}}^{-1}$ 
above and the similarity coordinate $\xi=x/t$ reads 
\begin{equation}
   \mathbf{Q}_{\textnormal{RP}}(\xi) = \frac{1}{2} \tilde{\mathbf{R}} \left( \mathbf{I} + \textnormal{sign}(\tilde{\boldsymbol{\Lambda}} - \mathbf{I} \xi) \right) \tilde{\mathbf{R}}^{-1} \cdot \mathbf{Q}_L + 
	                   \frac{1}{2} \tilde{\mathbf{R}} \left( \mathbf{I} - \textnormal{sign}(\tilde{\boldsymbol{\Lambda}} - \mathbf{I} \xi) \right) \tilde{\mathbf{R}}^{-1} \cdot \mathbf{Q}_R,    
										\label{exact.riemann} 
\end{equation} 
with $\mathbf{I}$ being the identity matrix. From $\mathbf{Q}_{\textnormal{RP}}(\xi)$ we can obtain the following 
\textit{Godunov state} $\mathbf{Q}_{\textnormal{God}} = \mathbf{Q}_{\textnormal{RP}}(0)$ at the interface ($\xi = 0$ ) 
\begin{eqnarray}
\mathbf{Q}_{\textnormal{God}} & = & \left(0
		,\frac{\sigma_{xx}^L c_p^2+2\sigma_{xx}^L c_s^2+ \sigma_{yy}c_p^2}{c_p^2}
		,\frac{\sigma_{xx}^L c_p^2+2\sigma_{xx}^L c_s^2+ \sigma_{zz}c_p^2}{c_p^2}
		,0
		,\sigma_{yz}^L  
		,0,\right. \nonumber \\ && \left.
		 \frac{c_p \rho u^L-\sigma_{xx}^L}{c_p \rho}
		,\frac{c_s \rho v^L-\sigma_{xy}^L}{c_s \rho}
		,\frac{c_s \rho w^L-\sigma_{xz}^L}{c_s \rho},\lambda, \mu, \rho, 1 \right), \nonumber
\end{eqnarray}
from which it is clear that all the components of the normal stress in $x$-direction ($\sigma_{xx},\sigma_{xy}$ and $\sigma_{xz}$) are zero, which means that the free surface 
boundary condition $\boldsymbol{\sigma} \cdot \mathbf{n}$ is indeed respected. 
 
As one can note, the model (\ref{eq:3.00})-(\ref{eq:3}) involves divisions by $\alpha$ that can be a source of instabilities at the interface, since the color function $\alpha$  is ideally set to zero or at least close to zero outside the solid medium. In order to address this problem, we introduce a simple transformation that avoids the divisions by zeros. In particular, we substitute all multiplications by $\alpha^{-1}=1/\alpha$, with  
\begin{eqnarray}
\alpha^{-1}\cong\frac{\alpha}{\alpha^2+\epsilon(\alpha)},
\label{eq:trick1}
\end{eqnarray}
where $\epsilon=\epsilon(\alpha)$ has to satisfy $\epsilon(1)=0$ and $\epsilon(0)=\epsilon_0>0$ in order to be consistent with the linear elasticity equations. In our case we take a simple linear function $\epsilon=\epsilon_0 (1-\alpha)$ with $\epsilon_0=10^{-3}$. 
The introduction of this new parameter with this method is mandatory to obtain a stable solution. 
The new eigenvalues are $\tilde{\lambda}=f \lambda$, where $f=\frac{\alpha}{\sqrt{\alpha^2+ \epsilon_0(1-\alpha)}}$ that for $\alpha \in [0,1]$ satisfies 
$f \in [0,1]$ and $f=1$ for $\alpha=1$.

As soon as we use a non-trivial geometry we obtain a diffuse interface of \textit{finite width} for the transition between the solid medium $\alpha=1$ and the surrounding gas / vacuum ($\alpha=0$).
For a relatively large width of the diffuse interface, there are some questions that arise naturally concerning the distribution of the characteristic function $\alpha$ inside the diffuse interface and the resulting effective position of the free surface boundary. 
In general, it is important to set up the diffuse interface shape such that $\nabla \alpha$ is oriented as the normal vector to 
the physical surface, i.e.~$\nabla \alpha \approx \mathbf{n}$. A simple way to do this is to represent the transition region by a piecewise polynomial. Let $r=r(\xx)$ be the signed distance between the real physical  
interface location and a generic point $\xx$ under consideration.  We then define the shape of the diffuse interface as function of a finite interface thickness $I_D \geq 0$, a shifting parameter $\eta$ and the auxiliary function:  
\begin{eqnarray}
\xi(r)=\left\{ 
\begin{array}{lll}
	1 & \textnormal{if} & r > \phantom{-} (1+\eta) I_D, \\ 
	0 & \textnormal{if} & r < -(1-\eta) I_D, \\
	\frac{r + (1-\eta) I_D}{2 I_D} & \textnormal{if} & r \in [-(1-\eta) I_D, (1+\eta) I_D]. 
\end{array}
\right.
\label{eq:xi}
\end{eqnarray}
We finally define the solid volume fraction as 
\begin{equation}
\alpha(r)=(1-\xi(r))^{p_d},
\label{eqn.alpha}
\end{equation} 
where  $p_d > 0$ is an exponent that determines the shape of the diffuse interface. 

\section{Numerical scheme}\label{sec:numerical-scheme}
The numerical method that we use in order to solve the PDE system introduced in Section \ref{sec_MatMod} is an explicit ADER-DG scheme that is of arbitrary high order accurate in space and time. The numerical method was presented for different PDE systems in \cite{DumbserPeshkov2016,Zanotti2015}, hence in the following we only give a brief summary. The PDE system \eqref{eq:3} can be written in compact matrix vector notation as 
\begin{eqnarray}
		\diff{\Q}{t}+\B(\Q)\cdot \nabla \Q = \S(\xx,t),
\label{eq:nm_0}
\end{eqnarray} 
where $\Q$ is the state vector, $\B(\Q) \cdot \nabla \Q$ is a non conservative product (see \cite{DLMtheory,Castro2006,Pares2006}) and $\S(\xx,t)$ is a known source term. In regions where $\alpha=1$ and thus $\nabla \alpha=0$, the PDE system  \eqref{eq:nm_0} reduces to the classical linear elastic wave equations \eqref{eq:1}, while for $\nabla \alpha \neq 0$ the system becomes locally \textit{nonlinear} and therefore requires a very robust numerical scheme as well as high resolution to be properly solved. Within this paper we use the simple and very robust  sub-cell finite-volume limiter approach in combination with  adaptive mesh refinement (AMR). A detailed description of the limiter can be found in \cite{Dumbser2014}.

In this section we report a short overview of the numerical scheme adopted in this paper, for more details we refer to \cite{DumbserPeshkov2016,Zanotti2015}. We discretize a $d$-dimensional computational domain $\Omega$ with a Cartesian grid as
\begin{eqnarray}
\Omega= \bigcup\limits_{i=1}^{N_e} T_i,
\label{eq:nm_1}
\end{eqnarray}
where $N_e$ is the total number of elements that in a Cartesian mesh can always be written dimension by dimension as $N_e=I_{max}\cdot J_{max} \cdot K_{max}$.
Since we are interested in a high order scheme, we first define a piecewise polynomial nodal basis $\{\phi_k\}_{k=1\ldots (N+1)^d}$ as the set of Lagrange polynomials passing through the Gauss-Legendre quadrature points on a reference unit element $T_{ref}$ for a given polynomial degree $N\geq0$. 
A weak formulation of the PDE system is obtained after multiplying Eq. \eqref{eq:nm_0} by a test function $\phi_k$ for $k=1\ldots (N+1)^d$ and then integrating over a space-time control volume $T_i \times [t^n, t^{n+1}]$: 
\begin{eqnarray}
\int\limits_{t^n}^{t^{n+1}}\int\limits_{T_i}\phi_k \left(\diff{\Q}{t}+\B(\Q)\cdot \nabla \Q \right) \dxt=\int\limits_{t^n}^{t^{n+1}}\int\limits_{T_i}\phi_k \S(\xx,t) \dxt.
\label{eq:nm_2}
\end{eqnarray}
We restrict the discrete solution to the space of piecewise polynomials of degree $N$, i.e.~the numerical solution $\u_h$ is written inside each element in terms of the polynomial basis as
\begin{eqnarray}
\u_h(\xx,t^n)=\sum\limits_{k=1}^{(N+1)^d}{\phi_k(x) \uh^n_k}=\bm{\phi}(\xx)\cdot \uh^n,
\label{eq:nm_3}
\end{eqnarray}
for $\xx \in T_i$ and $i=1 \ldots N_e$. The vector of degrees of freedom of $\u_h(\xx,t^n)$ is denoted by $\uh^n$. Using the definition \eqref{eq:nm_3} in the weak formulation given by Eq. \eqref{eq:nm_2} we obtain
\begin{eqnarray}
\left(\int\limits_{T_i}\phi_k \phi_l \dx \right) \left( \uh^{n+1}_l-\uh^n_l \right) + \int\limits_{t^n}^{t^{n+1}}\int\limits_{\partial T_i} \phi_k \mathcal{D}(\q_h^{-},\q_h^{+})\cdot \mathbf{n} \, dS \, dt + \int\limits_{t^n}^{t^{n+1}}\int\limits_{T_i} \phi_k \B(\q_h) \cdot \nabla \q_h \,  \dxt=\int\limits_{t^n}^{t^{n+1}}\int\limits_{T_i}\phi_k \S(\xx,t) \,  \dxt,
\label{eq:nm_4}
\end{eqnarray}
where we have introduced the jump contribution $\mathcal{D}(\q_h^{-},\q_h^{+})\cdot \mathbf{n} $ on the element boundaries and a the space-time  predictor solution $\q_h(\xx,t)$. More details concerning the computation of  $\q_h(\xx,t)$ will be reported later. 
For the approximation of the jump term $\mathcal{D}$ we use a path conservative scheme as introduced by Par\'es in \cite{Pares2006} and Castro et. al. in \cite{Castro2006}. We introduce a Lipschitz continuous path function  $\psi(\q_h^{-},\q_h^{+},s)$ defined for $s \in [0,1]$ such that $\psi(\q_h^{-},\q_h^{+},0)=\q_h^{-}$ and  $\psi(\q_h^{-},\q_h^{+},1)=\q_h^{+}$, where $\q_h^-$ denotes the boundary-extrapolated state from within the element $T_i$ and $\q_h^+$ the boundary-extrapolated state from the neighbor element. The simplest possible choice for $\psi$, which we use in this paper, is the linear segment path between the two states $\q_h^{-}$ and $\q_h^{+}$:
\begin{eqnarray}
		\psi(\q_h^-,\q_h^+,s) = \q_h^- + s \left(\q_h^+ - \q_h^- \right).
\label{eq:segmpath}
\end{eqnarray}
Following \cite{Pares2006,Castro2006} we now write the jump contribution as 
\begin{eqnarray}
\mathcal{D}(\q_h^{-},\q_h^{+})\cdot \mathbf{n}= \frac{1}{2}\left(\int\limits_0^1 \B(\psi(\q_h^{-},\q_h^{+},s))\cdot \mathbf{n} \, ds \right) \left(\q_h^+ - \q_h^- \right), 
\label{eq:nm_5}
\end{eqnarray}
which has to obey the general Rankine-Hugoniot condition \cite{DLMtheory,Pares2006}
\begin{eqnarray}
\mathcal{D}(\q_h^{-},\q_h^{+})\cdot \mathbf{n}-\mathcal{D}(\q_h^{+},\q_h^{-})\cdot \mathbf{n}=\int\limits_0^1 \B(\psi(\q_h^{-},\q_h^{+},s))\cdot \mathbf{n}\diff{\psi}{s} ds.
\label{eq:nm_6}
\end{eqnarray}
The previous integral can simply be evaluated \textit{numerically} using a sufficient number of Gaussian quadrature points. As Riemann solver we use the new HLLEM-type Riemann solver for non-conservative systems recently described in \cite{NCP_HLLEM}, since we want to preserve exactly the material parameters that appear in the PDE system.

Regarding the space-time predictor, we need to introduce a new polynomial basis of degree $N$ in space and time $\{\theta_k\}_{k=1\ldots (N+1)^{d+1}}$ where now $\theta_k(\xx,t) \in T_i \times [t^n,t^{n+1}]$ contains also the time. We represent $q_h(\xx,t)$ in terms of this new space-time basis as
\begin{eqnarray}
\q_h(\xx,t)=\sum\limits_{k=1}^{(N+1)^{d+1}}{\theta_k(\xx,t) \hat{\q}^n_k}.
\label{eq:nm_3_st}
\end{eqnarray}
Let $T_i^\circ=T_i-\partial T_i$ denote the interior of $T_i$ and $T_i^{st}=T_i^\circ \times [t^n,t^{n+1}]$ denote the new space-time control volume. The space-time predictor is then computed as an element-local solution of the following weak formulation of the PDE system \eqref{eq:nm_0}:
\begin{eqnarray}
\int\limits_{T_i^{st}}{\theta_k \diff{\q_h}{t} \dxt} + \int\limits_{T_i^{st}}{\theta_k \B(\q_h) \cdot \nabla \q_h \dxt} = \int\limits_{T_i^{st}}{\theta_k \S(\xx,t) \dxt} .
\label{eq:nm_7}
\end{eqnarray}  
%Since $\theta_k$ depends explicitly by the time, we 
%Integrating the first term of \eqref{eq:nm_7} by parts and using upwind on the lower slice $t^n$
% that allows to connect $q_h$ to the solution at the previous time step $u_h^n$:
Using integration by parts in the first term of Eq. \eqref{eq:nm_7} we obtain two spatial contributions on $T_i$ at $t^{n+1}$  and $t^n$ and an internal one since $\theta_k=\theta_k(\xx,t)$ contains explicitly the time. For the spatial contribution at time $t^n$ we use the numerical solution from the previous time step. Notice that this  corresponds to upwinding in the time direction due to the causality principle: 
\begin{eqnarray}
\int\limits_{T_i}{\theta_k(\xx,t^{n+1}) \q_h(\xx,t^{n+1})\dx}-\int\limits_{T_i}{\theta_k(\xx,t^{n}) \u_h(\xx,t^{n})\dx}-\int\limits_{T_i^{st}}{\diff{\theta_k(x,t)}{t} \q_h(\xx,t)\dxt} && \nonumber \\ +\int\limits_{T_i^{st}}{\theta_k \B(\q_h) \cdot \nabla \q_h \dxt} = \int\limits_{T_i^{st}}{\theta_k \S(\xx,t) \dxt}. && 
\label{eq:nm_8}
\end{eqnarray}
Since Eq. \eqref{eq:nm_8} is element-local it can be solved using a simple and efficient Picard method without any communication with the neighbor elements, see e.g. Dumbser et. al. \cite{Dumbser2008}.

The numerical scheme is constrained by a \textit{local} CFL-type stability condition, see \cite{Dumbser2008,AMR3DCL,Zanotti2015}, that is given by  
\begin{eqnarray}
\Delta t < \frac{\textnormal{CFL}}{d}\frac{h}{2N+1}\frac{1}{|\lambda_{max}|}, 
\label{eq:cfl}
\end{eqnarray}
where $h$ is the local mesh size, $\lambda_{max}$ is the maximum eigenvalue of the PDE system, and CFL$<1$ is the  Courant number, which should be chosen according to \cite{Dumbser2008} in order to have linear stability. 
Concerning the adaptive mesh refinement (AMR) we rely on the ExaHyPE engine \textcolor{blue}{http://exahype.eu}, which is built in turn upon the space-tree implementation Peano \cite{Peano1,Peano2} 
realising cell-by-cell refinement \cite{Khokhlov1998}. For further details about AMR in combination with high order finite volume and DG schemes with time-accurate local time stepping (LTS), see 
\cite{AMR3DCL,Zanotti2015c,Zanotti2015d,ADERDGVisc}. 

In order to decide where to refine, we introduce a simple refinement indicator function named $\varphi=\varphi(\xx,t)$ that defines the observed variable for the refinement/recoarsening process and a so called real-valued {\itshape estimator function} $\chi=\chi[\varphi]$, see again \cite{AMR3DCL} for more details. After  defining the indicator function, we define the cell-averages of $\varphi$ as 
\begin{equation}
\hat{\varphi}_i = \frac{1}{|T_i|} \int_{T_i} \varphi (\xx,t)\,d\xx \qquad \qquad \forall i=1\ldots N_e,
\end{equation}  
and then we compute the estimator function as
\begin{equation}
\chi_i[\varphi] = \max_{c \in \mathcal{V}_i} \left( \left | \hat{\varphi}_c  -\hat{\varphi}_i \right|  /  \left\| \mathbf{x}_c - \mathbf{x}_i \right\| \right),
\end{equation}
where $\mathcal{V}_i$ contains all the Voronoi elements of $i$. Our estimator function $\chi$ is simply based on 
an approximation of the gradient of the solution in several spatial directions \cite{AMR3DCL}. With these ingredients at hand, we introduce a simple rule for the refinement/recoarsening process based on two thresholds $\chi^+$ and $\chi^-$ as follows: 
\begin{enumerate}[(i)]
	\item if $\chi_i[\varphi] > \chi^+$ then $T_i$ is labeled for mesh refinement;
	\item if $\chi_[\varphi] < \chi^-$ then $T_i$ is labeled for mesh recoarsening.
\end{enumerate} 
Within this paper, we always use $\varphi(\xx,t)=\varphi(Q)=\alpha$, $\chi^+=0.01$ and $\chi^+=0.001$. We will also use the volume fraction $\alpha$ to specify the zones where to activate the subcell finite volume limiter \cite{Dumbser2014}. In particular, we activate the FV limiter whenever $ \alpha \notin [\epsilon,1-\epsilon]$, 
with $\epsilon=10^{-3}$.
Since the topology of the geometry described by $\alpha$ is supposed to be stationary in time, we can consider the refinement and the limited zones also as steady and therefore they need to be identified only once in the mesh initialization step.

\section{Numerical results}\label{sec:numerical-results}

\subsection{Reflected plane wave}
\begin{figure}[ht!]%
\begin{center}
\includegraphics[width=0.6\columnwidth]{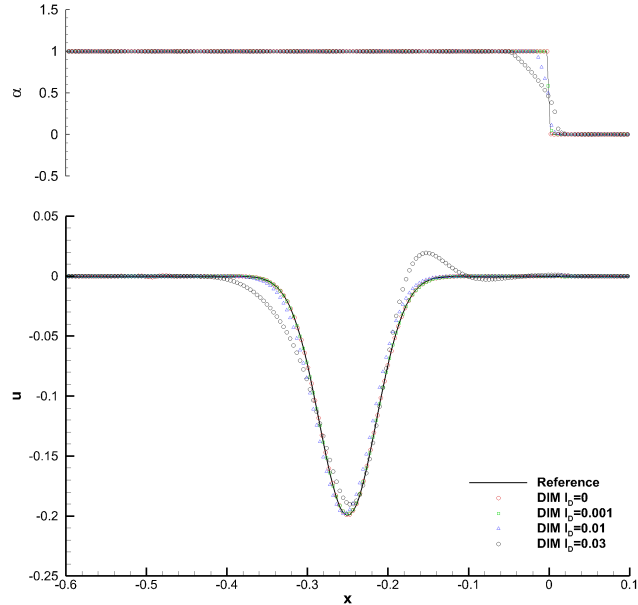} 
\end{center}
\caption{ (left) Numerical results obtained with the new diffuse interface approach for a plane wave reflection  problem on a free surface located in $x=0$ using a variable interface thickness of $I_D=0$, $I_D=0.001$, $I_D=0.01$ and $I_D=0.03$. In all four cases we report the velocity component $u$ compared with the exact solution of the problem (bottom) together with the spatial distribution of $\alpha$ (top). } 
\label{fig.2}%
\end{figure}

The purpose of this first test problem is to systematically study the influence of the width $I_D$ of the diffuse interface layer onto the numerical results. We also 
show that the model indeed converges to the correct solution in the limit $I_D \to 0$. 
We take a simple plane wave impulse in a domain $\Omega=[-1,1]\times[-0.1,0.1]$ initially placed at  $x_0=-0.25$ and hitting a free surface boundary placed in $x_D=0$. The Lam\'e constants are chosen as  
$\lambda=2$, $\mu=1$ and $\rho=1$. We define $\mathbf{Q}_0=( 0, 0, 0, 0, 0, 0, 0, 0, 0, \lambda, \mu, \rho, \alpha(x))$ and $\boldsymbol{\delta}=(0.4, 0.2, 0.2, 0, 0, 0, -0.2, 0, 0,  0, 0, 0, 0 )$ and set 
$$\mathbf{Q}(x,y,t=0)= \mathbf{Q}_0 + \boldsymbol{\delta} \cdot e^{-\frac{(x-x_0)^2}{\epsilon^2}},$$ 
with the halfwidth $\epsilon=0.05$. The volume fraction function $\alpha(x)$ is prescribed according to \eqref{eqn.alpha} and \eqref{eq:xi}. 
We use an ADER-DG P$_4$ scheme and a uniform Cartesian grid with  $100 \times 2$ elements. The mesh resolution is chosen fine enough so that the numerical results are grid-independent and only
depend on the choice of the interface thickness $I_D$.  
Since for this test $c_p=2$, the exact solution at time $t=t_{end}=0.25$ is the reflected p-wave which is located again in the initial position. 
We consider four cases with different choices of the interface width $I_D$, ranging from $I_D=0.03$ to the limit $I_D=0$,
where the interface is exactly located on a cell boundary. 
From the results depicted in Figure \ref{fig.2} we can conclude that the diffuse interface method is able to 
reproduce the exact solution of the problem for sufficiently small values of the interface thickness $I_D$. 
We also stress that the use of a path-conservative method allows us to reduce the interface thickness 
exactly to $I_D=0$, which leads to a jump in $\alpha$ at an element interface, but which is still 
properly accounted for thanks to the jump terms $\mathcal{D}_{i \pm \frac{1}{2}}$ used in the numerical scheme. 

For rather large values of the finite interface thickness $I_D$, where the actual shape of the spatial distribution of $\alpha$ starts to play a role, we have found empirically that a good choice for the parameters $\eta$ and $p_d$  in \eqref{eq:xi} is $\eta=-0.6$ and $p_d=0.5$. This choice allows to obtain still a correct reflection of a p-wave  even for very thick interfaces. However, for sufficiently small values of $I_D$, the choice of $\eta$ and $p_d$ has only very little influence. 

\subsection{Scattering of a plane wave on a circular cavity}\label{sec:scattering-of-a-plane-wave-on-a-circular-cavity}
\begin{figure}[ht!]
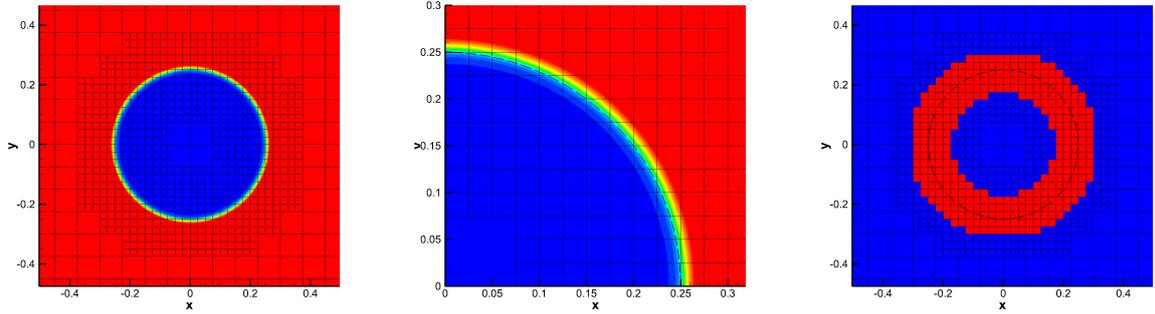
%
	\begin{center}
		\includegraphics[width=0.32\columnwidth]{\flink{}{PW_alpha}} 
		\includegraphics[width=0.32\columnwidth]{\flink{}{PW_alpha_2}} 
		\includegraphics[width=0.32\columnwidth]{\flink{}{PW_alpha_limiter}}
	\end{center}
	\caption{Setup of the scattering test problem. AMR grid and distribution of the characteristic function $\alpha$ (left). Detail of the free surface location $\partial C$ shown via a dashed line and $\alpha$ color contours (center). Limited cells highlighted in red and unlimited cells shown in blue (right). }%
	\label{fig.pw1}%
\end{figure}

In this test case we consider an initially planar $p$-wave traveling in $x$-direction inside a solid medium and hitting an empty circular cavity. The initial state is given by 
\begin{eqnarray}
	\mathbf{Q}(x,0) = (0,0,0,0,0,0,0,0,0, \lambda, \mu, \rho, \alpha) + 0.1\cdot(4,2,2,0,0,0,-2,0,0,0,0,0,0)\sin(2 \pi x),
\label{eq:pw1}
\end{eqnarray}
with $\lambda=2$, $\mu=1$ and $\rho=1$. 
The value of $\alpha$ is parameterised through the circular surface $C=\{(x,y)\,\, | \,\, x^2+y^2  \leq 0.25^2\}$ so that $\alpha(\mathbf{x})=0$ if $\mathbf{x} \in C$ and $\alpha=1$ if 
$\mathbf{x} \notin C$. The width parameter of the diffuse interface is set to $I_D=0.01$ on $\partial C$. % with respect to the direction $\vec{n}=\vec{r}$. 
The computational domain is $\Omega=[-3,3]^2$ and the initial Cartesian grid at level $\ell=0$ consists of  $80 \times 80$ cells.  We then use one further refinement level $\ell_{max}=1$ based on the gradient of $\alpha$ in order to refine the  mesh close to the diffuse interface. Furthermore we use a fifth order ADER-DG method based on piecewise polynomials of degree $N=4$ in both space and time, supplemented with a second order TVD sub-cell finite volume limiter. The resulting AMR grid and the color contours of $\alpha$ are shown in Figure \ref{fig.pw1}, together with the region where the subcell finite volume limiter is activated. From the plot in the central panel of Figure \ref{fig.pw1}
one can see that the width of the interface layer is of the order of the size of one cell of the high order 
DG scheme. 

\begin{figure}[ht!]
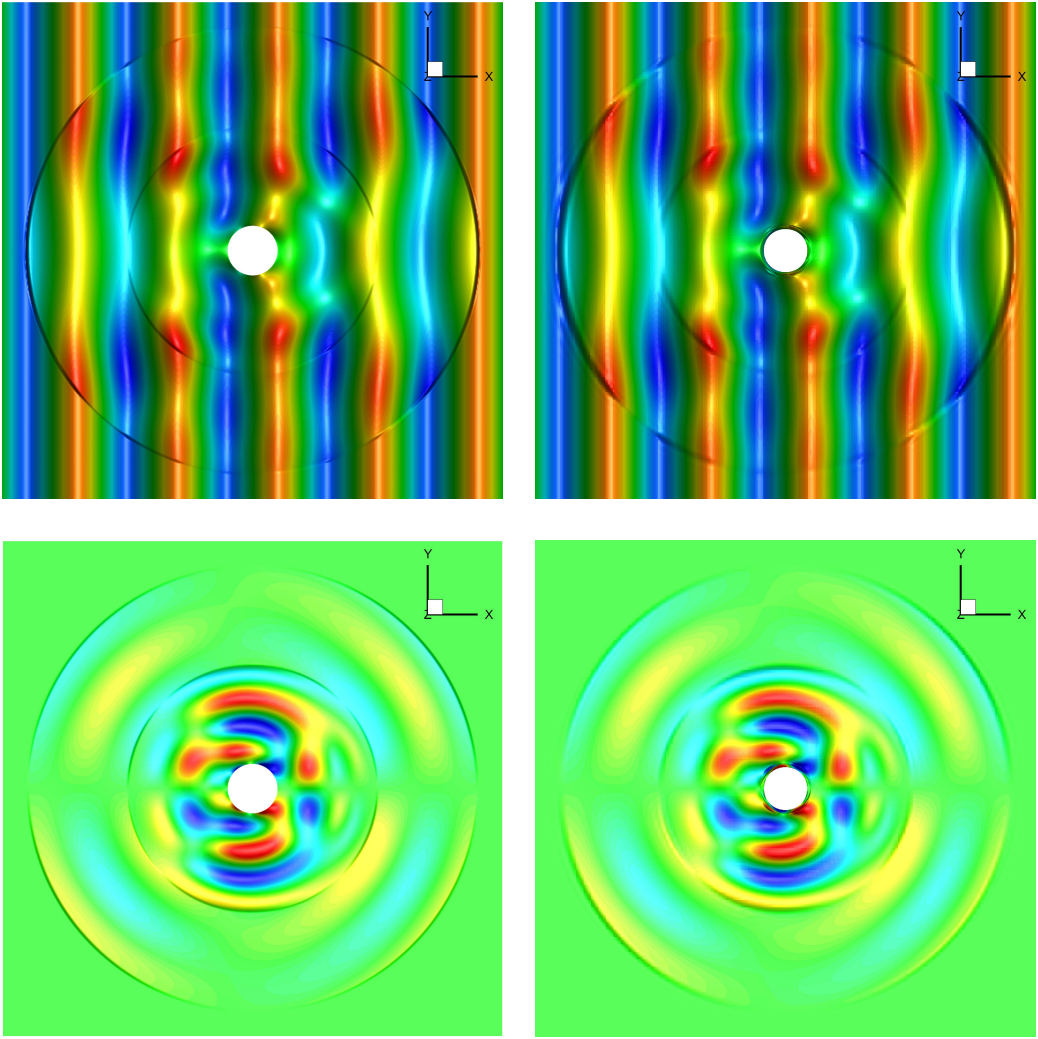
%
	\begin{center}     		
		\begin{tabular}{lr} 
		\includegraphics[width=0.4\columnwidth]{\flink{}{PW_sxx_PDESol}} & 
		\includegraphics[width=0.4\columnwidth]{\flink{}{PW_sxx_DIM}}    \\ \\ 
		\includegraphics[width=0.4\columnwidth]{\flink{}{PW_sxy_PDESol}}  & 
		\includegraphics[width=0.4\columnwidth]{\flink{}{PW_sxy_DIM}} 
		\end{tabular} 
	\end{center}
	\caption{2D wave scattering problem at time $t=1$. Reference solution (left) and solution obtained with the new  diffuse interface approach on AMR grid (right) for $\sigma_{xx}$ (top) and $\sigma_{xy}$ (bottom).}  
	\label{fig.pw3}%
\end{figure}

\begin{figure}[ht!]
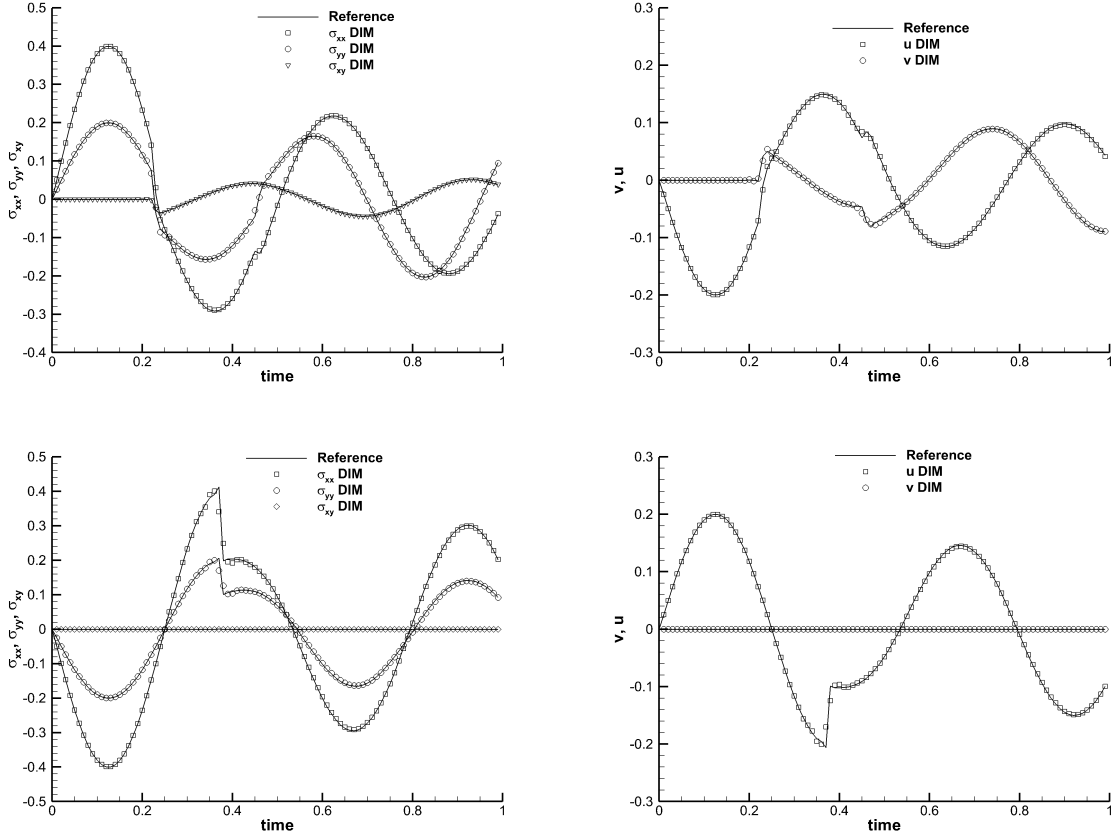
%
	\begin{center}
		\includegraphics[width=0.48\columnwidth]{\flink{}{PW_TS_1s}}
		\includegraphics[width=0.48\columnwidth]{\flink{}{PW_TS_1v}} \\
		\includegraphics[width=0.48\columnwidth]{\flink{}{PW_TS_2s}}
		\includegraphics[width=0.48\columnwidth]{\flink{}{PW_TS_2v}} \\
	\end{center}
	\caption{2D wave scattering problem. Comparison between the reference solution and the numerical results  obtained with the new diffuse interface method  on AMR grid for two seismograms recorded in $\mathbf{x}_1=(0.5,0.5)$ (top row) and $\mathbf{x}_2=(1.0,0.0)$ (bottom row)}%
	\label{fig.pw4}%
\end{figure}

In Figure \ref{fig.pw3} we report the numerical results obtained with the new diffuse interface approach on adaptive Cartesian grids at time $t=1$ and compare them with a reference solution that has been obtained with a third order ADER-DG scheme ($N=2$) on a very fine unstructured boundary-fitted mesh \cite{gji1} composed of  $N_e=563,280$ triangles. 
Figure \ref{fig.pw4} shows a comparison between the reference solution and the numerical solution obtained with the diffuse interface method via numerical seismograms that have been recorded in two receiver locations $\xx_1=(0.5,0.5)$ and $\xx_2=(1.0,0.0)$. A very good agreement between the new diffuse interface method and the reference solution is obtained for this test case. 

At this point we would like to stress again that in the new diffuse interface approach, the presence of the boundary  condition is included in the PDE system \textit{only} by choosing a spatially variable value 
of $\alpha$. The AMR grid is \textit{not}   
at all aligned with the free surface boundary and remains always locally Cartesian (with $h$-adaptivity).  Furthermore, the time step size in our approach is \textit{not} affected by the so-called small cell problem or sliver element problem, as it would have been the case for Cartesian cut-cell methods or low quality unstructured meshes and which usually requires a special treatment \cite{gij5,SIStag2017}. In our diffuse  interface approach, the eigenvalues of the PDE system are independent of $\alpha$ and also our mesh 
can be chosen independently of $\alpha$ and almost independently of the geometry of the problem to be solved (apart from local $h$ adaptivity used in regions of strong gradients of $\alpha$). Therefore, the admissible local time step size is only governed by the maximum wave  speed $c_p$ and the local mesh size of the AMR grid, and not by the geometry of the problem to be solved. Note that in all our simulations on AMR grids, we use time-accurate  local time stepping (LTS), see  \cite{gij5,AMR3DCL,AMR3DNC,Zanotti2015} for details. 

\subsection{2D tilted Lamb problem}
\label{NS:2DLamb}
In this test case we to study the two dimensional tilted Lamb problem. We take the same setup as used in \cite{Komatitsch1998,gji1,SIStag2017}. The physical domain $\Omega=\{(x,y)\in \R^2 \, | \, 0 \leq x\leq 4000 \,\, , \,\, 0 \leq y\leq 2000+\tan{(\theta)}x\}$ contains a free surface with a tilt angle of $\theta=10^\circ$, so that the boundary is \textit{not} grid aligned when using a Cartesian mesh along the coordinate axes. The  computational domain used for the diffuse interface approach, however, is a simple rectangular box that fully contains $\Omega$. The initial Cartesian grid on the coarsest level $\ell=0$ has $ 96 \times 90$ cells and we use an ADER-DG $P_3$ scheme with subcell finite volume limiter to solve this problem. The chosen $p-$ and  $s-$wave velocities are set to $c_p=3200m\,s^{-1}$ and 
$c_s=1847.5m\,s^{-1}$, respectively. The mass density is taken as $\rho=2200 kg \, m^{-3}$ so that the resulting Lam\'e constants are $\lambda=7.5096725\cdot 10^{9}$ and $\mu=7.50916375\cdot 10^{9}$. 
The initial condition is $\mathbf{Q}(\mathbf{x},0)=0$ everywhere in $\Omega$. The wave propagation is driven by a directional point source located in $\xx_s=(1720.0,2265.28)$. We place two receivers, one close to the 
interface but slightly below, so that $\alpha=1$, $\mathbf{x}_2=(2694.96, 2460.08)$ and the other one exactly at 
the physical interface in $\mathbf{x}_1=(2694.96, 2475.08)$. As reference solution we use again an ADER-DG method on boundary-fitted unstructured meshes, which has already been carefully validated against the exact solution of Lamb's problem 
in \cite{gji1}. The reference solution is computed using a polynomial approximation degree $N=4$ in space and time and an unstructured mesh of $N_e=844,560$ triangles. The point source   
$$S(\xx,t)=\frac{1}{\rho}\vec{\bm{d}}\delta(\xx - \xx_s)\mathcal{S}(t)$$ 
is a delta distribution in space located in $\xx=\xx_s$ and its temporal part is a Ricker wavelet given by 
\begin{eqnarray}
	\mathcal{S}(t)=a_1\left( 0.5+a_2(t-t_D)^2 \right),
\label{eq:}
\end{eqnarray}
where $t_D=0.08s$ is the source delay time; $a_1=-2000kg\, m^{-2}\, s^{-2}$; $a_2=-(\pi f_c)^2$; and $f_c=14.5 Hz$.
Finally the vector $\vec{\bm{d}}=(-\sin\theta, \cos\theta ,0,0,0,0,0,0,0,0)^\top$ determines the direction of the impulse and takes into account the tilt angle $\theta$.
For this test we use an interface thickness of $I_D=2 m$. Furthermore we compare two different resolutions at the interface corresponding to a maximum refinement level of $\ell_{max}=2$ and $\ell_{max}=3$. 
%Furthermore, we take $\eta=0$ in order to take a look on how the shear waves, that dominates this problem, act with respect to a large interface. 
\begin{figure}[!h] 
	\begin{center}
		\includegraphics[width=0.8\columnwidth]{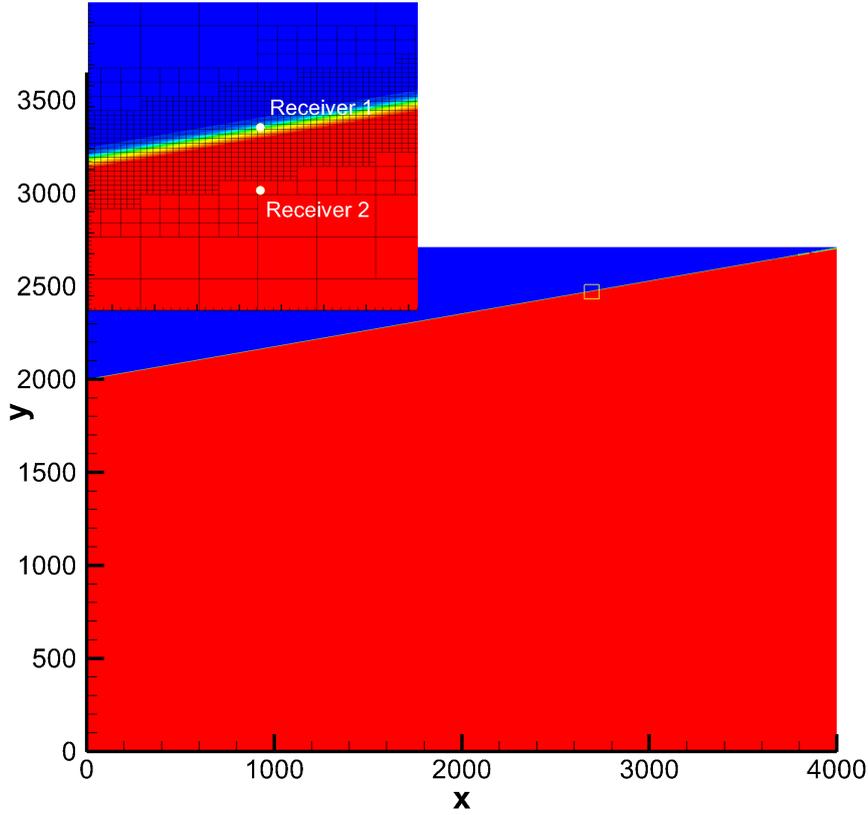} 
	\end{center}
	\caption{Tilted Lamb problem. Distribution of $\alpha$ in the computational domain and position of the two receivers. Note that the tilted free surface is not aligned with the Cartesian grid. 
	The resolution of the free surface is improved by a combination of AMR and subgrid finite volume limiter.  }%
	\label{fig.l2d1}%
\end{figure}
\begin{figure}[!h]
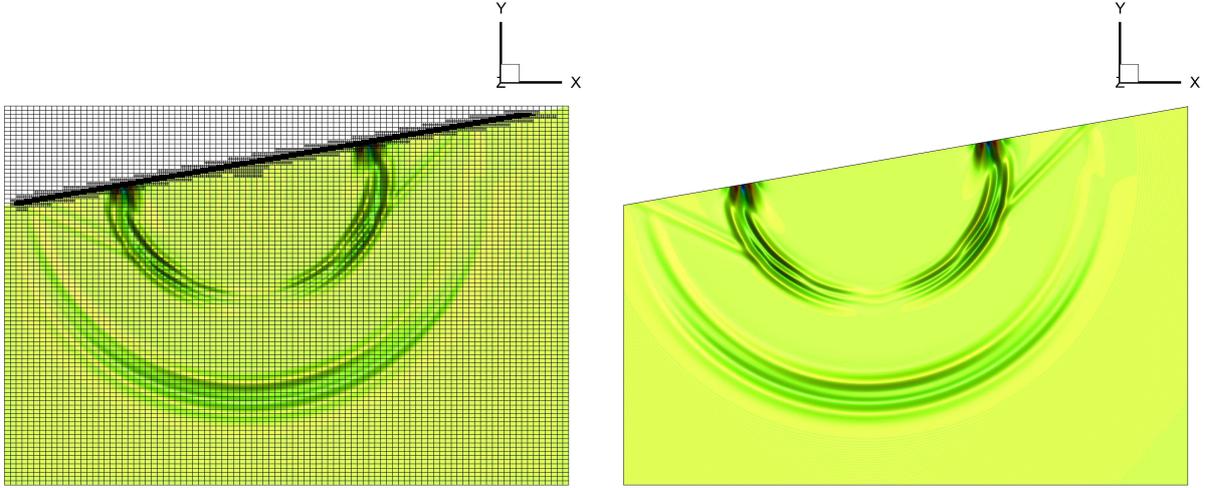
%
	\begin{center}
		\includegraphics[width=0.49\columnwidth]{\flink{}{L2D_vAMRng}}
		\includegraphics[width=0.49\columnwidth]{\flink{}{L2D_vPDESol}}
	\end{center}
	\caption{Tilted Lamb problem. Comparison of the vertical velocity $v$ between the new diffuse interface approach on AMR grid (left) and the reference solution obtained on a boundary-fitted onstructured mesh (right) at $t=0.6$}%
	\label{fig.l2d2}%
\end{figure}
\begin{figure}[!h]
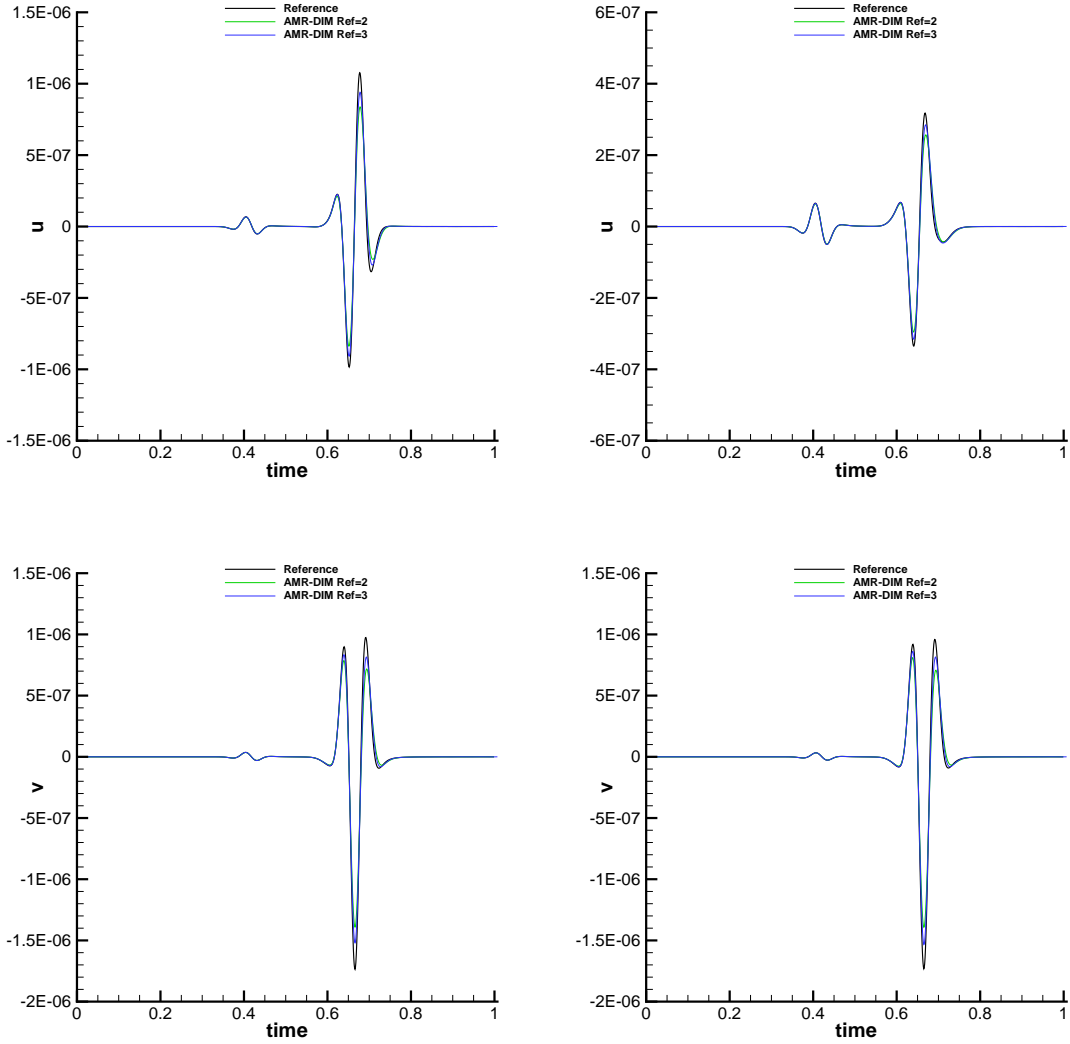
 
	\begin{center}
		\includegraphics[width=0.45\columnwidth]{\flink{}{L2D_u_pp1}}
		\includegraphics[width=0.45\columnwidth]{\flink{}{L2D_u_pp2}} \\
		\includegraphics[width=0.45\columnwidth]{\flink{}{L2D_v_pp1}}
		\includegraphics[width=0.45\columnwidth]{\flink{}{L2D_v_pp2}}
	\end{center}
	\caption{Lamb's problem. Comparison of the reference solution (solid black line) and the numerical solution obtained with the new diffuse interface approach on adaptive Cartesian grids for $2$ and $3$ 
	refinement levels in the two receivers 1 (left) and 2 (right).}%
	\label{fig.l2d3}%
\end{figure}
Figure \ref{fig.l2d1} shows the value of the solid volume fraction $\alpha$, as well as the positions of the seismogram recorders in $\mathbf{x}_1$ and $\mathbf{x}_2$. In Figure \ref{fig.l2d2} we compare between the 
numerical solution obtained with the new diffuse interface approach on Cartesian AMR grids and the reference solution obtained on a boundary-fitted unstructured mesh. We can observe a good agreement between the two 
solutions, which becomes also clear if we compare the seismograms, see Figure \ref{fig.l2d3}. In this case it 
is  also evident that the use of a higher grid resolution at the interface allows to approach the reference 
solution better. 

\subsection{Wave propagation in a complex 2D geometry}
\begin{figure}[!h]%
	\begin{center}\
		\includegraphics[width=0.8\columnwidth]{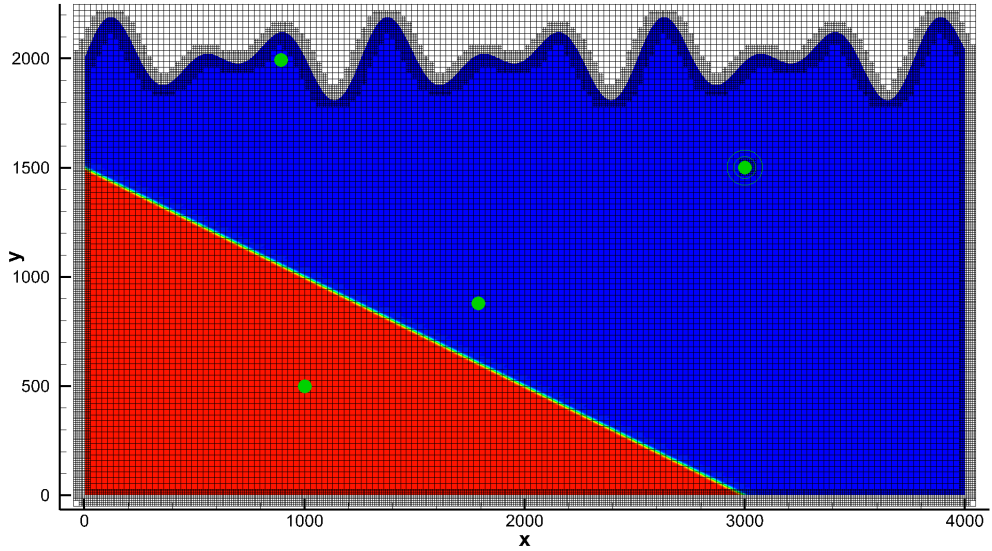}
	\end{center}
	\caption{Wave propagation in complex 2D geometry. Computational domain and AMR grid for the diffuse interface approach, colored by the mass density.}%
	\label{fig.CG1}%
\end{figure}

\begin{figure}[!h]
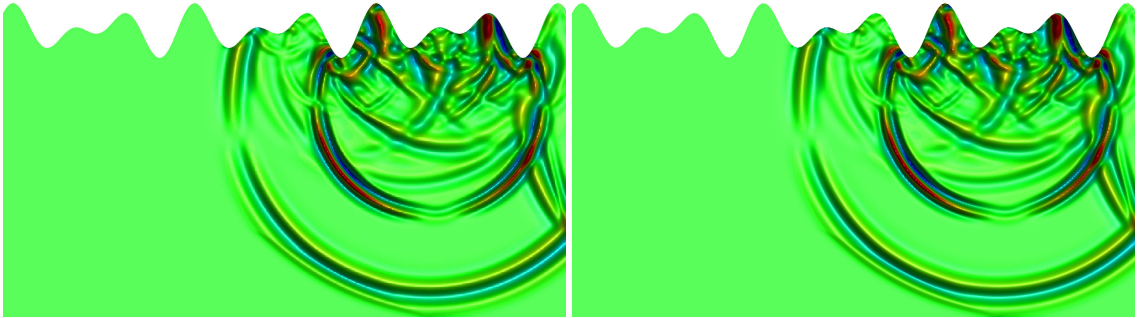
%
	\begin{center}
		\includegraphics[width=0.45\columnwidth]{\flink{}{CG_t4_AMR}}
		\includegraphics[width=0.45\columnwidth]{\flink{}{CG_t4_Ref}}
	\end{center}
	\caption{Wave propagation in complex 2D geometry. Comparison between the numerical solution obtained with the diffuse interface approach (left) and the reference solution (right) for $\sigma_{xx}$ at $t=0.5$. 
	         For the diffuse interface results, only the physically relevant part of the domain is shown. }%
	\label{fig.CG2}%
\end{figure}

This test is very similar to the tilted Lamb problem but in a non-trivial domain with heterogeneous material. The domain is $\Omega_f=\{(x,y) \,\, | \,\, x \in [0,4000] \,\, y \in [0,f(x)]\}$ where $f(x)=2000+100\left( \sin{(\frac{3}{200}x)}+ \sin{(\frac{2}{200}x)} \right)$ describes the upper topography. Free surface 
boundary conditions are imposed everywhere on $\partial \Omega_f$. 
\begin{table}[!t] 
\begin{center}
\begin{tabular}{cccc}
	Zone & $c_p$ (m/s) 	& $c_s$ (m/s) & Location	 \\
	\hline
	1		&	3200.00	& 1847.50 & $y>1500-\frac{x}{2}$ \\
	2		&	2262.74	& 1306.38 & $y\leq1500-\frac{x}{2}$ \\
	\hline
\end{tabular}
\end{center}
\caption{p- and s-wave speeds in the two layers used for the wave propagation problem in complex geometry.} 
\label{tab:LAYERS}
\end{table}
The heterogeneous material consists in two layers whose parameters are reported in Table \ref{tab:LAYERS}. The initial state vector is $\mathbf{Q}(x,0)=0$ and the wave propagation is driven by a point source placed in $\bar{\xx}=(3000,1500.18)$ as described in the previous Section \ref{NS:2DLamb}. Three seismometers are placed in the locations reported in Table \ref{tab:SYSCGEOM} and graphically depicted in Figure \ref{fig.CG1} to record the time history of the wave propagation. 
\begin{table}[!t]
\begin{center}
\begin{tabular}{cccc}
	Receiver & 1 	& 2  & 3 	 \\
	\hline
	$x$		&	893.80	& 1790.0 &  1000.0 \\
	$y$		&	1994.83	& 880.0  &  500.0  \\
	\hline
\end{tabular}
\end{center}
\caption{Receiver locations used for the seismogram recordings in the wave propagation problem in complex geometry.}
\label{tab:SYSCGEOM}
\end{table}
We take an extended domain $\Omega=[-50,4050]\times [-50, 2300]$ that fully contains $\Omega_f$. The initial Cartesian grid on level $\ell=0$ is composed of $ 160 \times 90$ elements. Subsequently, one refinement level is added in regions with large gradients of $\alpha$, i.e.~we set $\ell_{max}=1$. The value of $\alpha$ is used to  define the complex physical domain $\Omega_f$ following our diffuse interface approach. The chosen smoothing  parameter close to the upper surface is taken as $I_D=5.0$, based on the distance function from a point and the  boundary of the domain $\partial \Omega_f$. We furthermore use $I_d=0$ on the left, right and bottom boundaries,  which are all grid aligned. The resulting AMR grid and the spatial distribution of $\alpha$ in the computational domain are shown in Figure \ref{fig.CG1}. 
A direct comparison between the solution obtained with the novel diffuse interface approach using an ADER-DG $P_4$ scheme on the AMR grid and the reference  solution obtained with an ADER-DG  $P_4$ scheme on a boundary-fitted unstructured mesh composed of $20254$ triangles is reported in Figure \ref{fig.CG2}. 
The comparison of the seismograms at the three receivers up to $t=2.0s$ is shown in Figures \ref{fig.CG3u} and \ref{fig.CG3v}. A very good agreement is achieved for short times, and even at later times the agreement remains rather good, considering that at later times the signal is the result of several reflected waves on the free 
surface.
\begin{figure}[!h]%
	\begin{center}
		\includegraphics[width=0.9\columnwidth]{\flink{}{CG_ts_u_1}} \\ 
		\includegraphics[width=0.9\columnwidth]{\flink{}{CG_ts_u_2}} \\ 
		\includegraphics[width=0.9\columnwidth]{\flink{}{CG_ts_u_3}}   
	\end{center}
	\caption{Wave propagation in complex 2D geometry. Comparison of the velocity component $u$ at the three receivers $1-3$ from top to bottom. }% 
	\label{fig.CG3u}%
\end{figure}
\begin{figure}[!h]%
	\begin{center}
		\includegraphics[width=0.9\columnwidth]{\flink{}{CG_ts_v_1}} \\
		\includegraphics[width=0.9\columnwidth]{\flink{}{CG_ts_v_2}} \\
		\includegraphics[width=0.9\columnwidth]{\flink{}{CG_ts_v_3}}
	\end{center}
	\caption{Wave propagation in complex 2D geometry. Comparison of the velocity component $v$ at the three receivers $1-3$ from top to bottom. }% 
	\label{fig.CG3v}%
\end{figure}

\subsection{Scattering of a planar wave on a sphere} 
\begin{figure*}[!h]
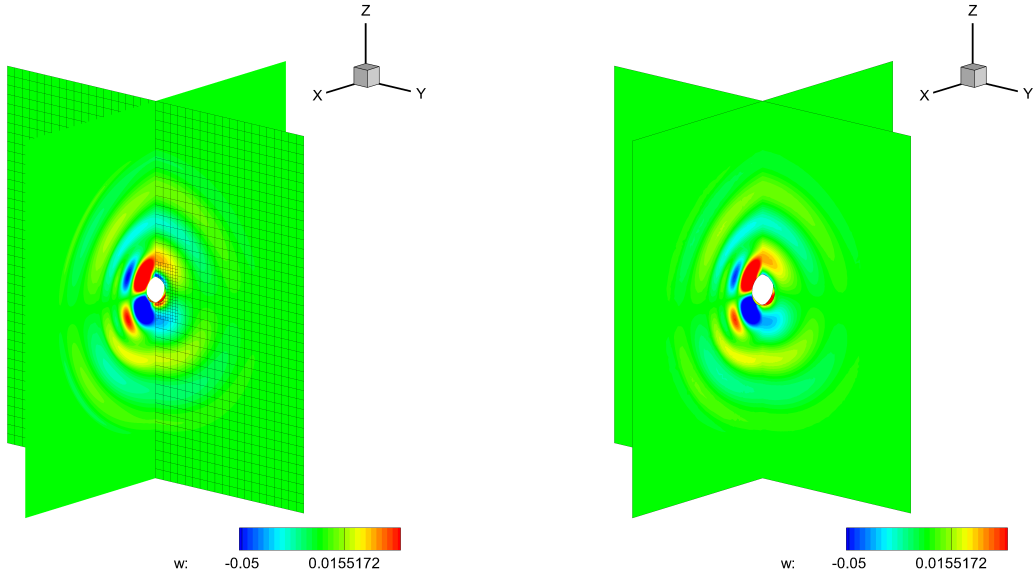

\includegraphics[width=0.48\columnwidth]{\flink{}{PW3D_AMRend}} 
\includegraphics[width=0.48\columnwidth]{\flink{}{PW3D_PDESOLend}} 
\caption{Scattering of a plane wave on a sphere. Velocity component $w$ at $t_{end}=1.0$ obtained with the new diffuse interface method on AMR grid (left) and the unstructured reference code \texttt{SeisSol} (right).}%
\label{fig.NT3DPW_2}%
\end{figure*}

\begin{figure*}[!h]
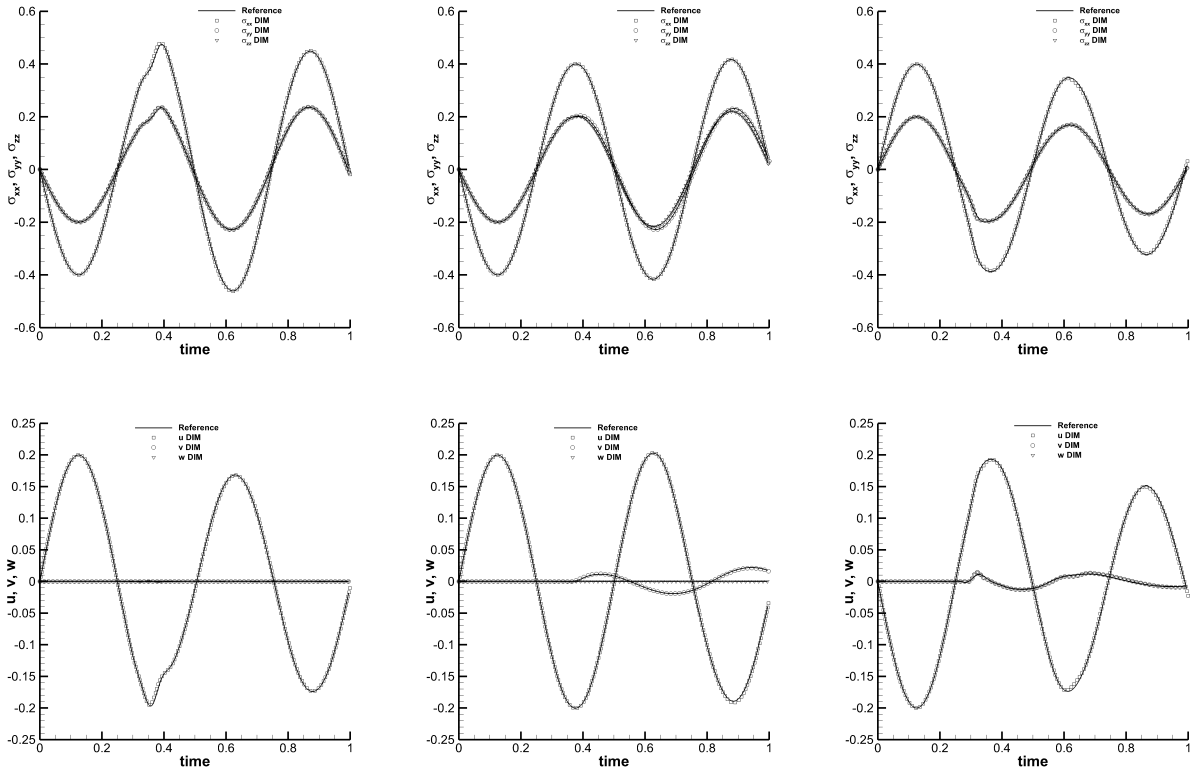

\includegraphics[width=0.33\columnwidth]{\flink{}{PW3D_ts_rec1_sigma}} 
\includegraphics[width=0.33\columnwidth]{\flink{}{PW3D_ts_rec2_sigma}} 
\includegraphics[width=0.33\columnwidth]{\flink{}{PW3D_ts_rec3_sigma}} \\
\includegraphics[width=0.33\columnwidth]{\flink{}{PW3D_ts_rec1_v}} 
\includegraphics[width=0.33\columnwidth]{\flink{}{PW3D_ts_rec2_v}} 
\includegraphics[width=0.33\columnwidth]{\flink{}{PW3D_ts_rec3_v}} 
\caption{Scattering of a plane wave on a sphere. Comparison of the resulting signal in the three receivers. In the first row we report the time series of the stress tensor components  
$\sigma_{xx}$, $\sigma_{yy}$ and $\sigma_{zz}$ for the receivers $1,2,3$, respectively, from left to right. In the second row the velocity signal is reported for the same receivers. }%
\label{fig.NT3DPW_1}%
\end{figure*}

Here we consider the 3D extension of the test reported in Section \ref{sec:scattering-of-a-plane-wave-on-a-circular-cavity} , which consists of a planar $p-$wave traveling in the $x-$direction and hitting a spherical cavity on which free surface boundary conditions apply. 
Our computational domain is the simple cube $\Omega=[-3,3]^3$ and the presence of the spherical obstacle is only taken into account by a spatially variable distribution of the volume fraction function $\alpha$. So $\alpha=1$ if $\alpha \notin B$ and $\alpha=0$ if $\alpha \in B$ where $B=\{(x,y,z) \,\, | \,\, x^2+y^2+z^2\leq 0.25^2\}$ is the sphere with radius $R=0.25$. The chosen interface width is $I_D=10^{-2}$. The computational domain is covered with a uniform initial mesh of $ 40 \times 40 \times 40$ elements. We then add one refinement level $\ell_{max}=1$ based on the gradient of $\alpha$. Furthermore we use piecewise polynomials of degree $N=5$ in space and time for this simulation. We consider three receivers placed in $\xx_1=(-1,0,0)$, $\xx_2=(0,-1,0)$ and $\xx_3=(0.5,0.5,0.5)$. 
As a reference solution we use again the explicit ADER-DG scheme implemented in the \texttt{SeisSol} code \cite{gij2,SeisSol1,SeisSol2} using a boundary-fitted unstructured grid with $N_e=31,732$ tetrahedral elements and piecewise polynomials of degree $N=4$ in space and time. \texttt{SeisSol} is a mature 
production code for large-scale seismic wave propagation problems in complex 3D geometries and has been 
heavily optimized so that it achieves a sustained Petaflop performance on modern supercomputers, see 
\cite{SeisSol1,SeisSol2} and \texttt{www.seissol.org}. 
A comparison of the contour colors for the velocity component $w$ is shown in Figure \ref{fig.NT3DPW_2} and a
direct comparison of the time series recorded in the three receivers is presented in Figure \ref{fig.NT3DPW_1}.  A very good agreement between the reference solution and the novel diffuse interface approach can be observed also in this case. 

\subsection{Wave propagation in a complex 3D geometry}
\begin{figure*}%
	\centering
	\includegraphics[width=0.7\columnwidth]{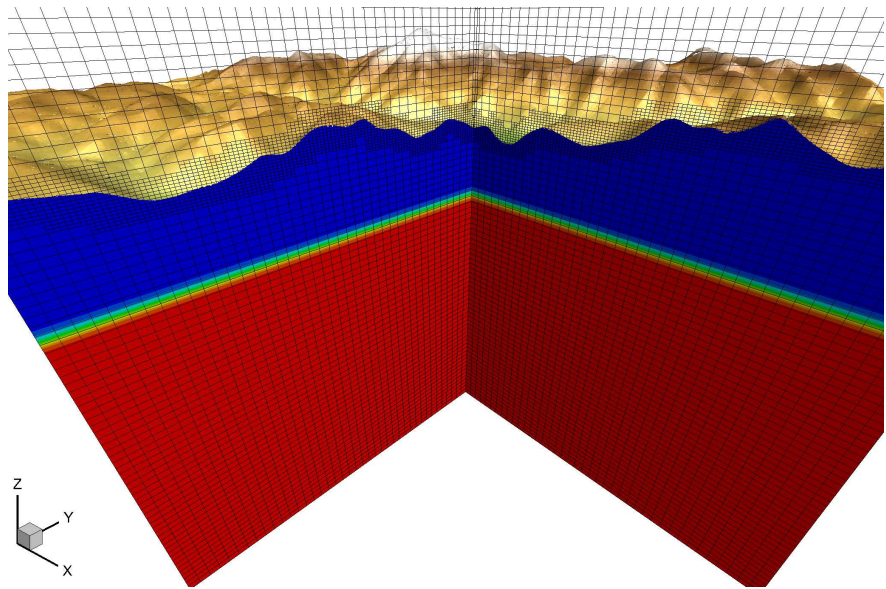} 
	\caption{Wave propagation in a realistic 3D geometry. Plot of the adaptive Cartesian mesh used for the test colored with the Lam\'e constant $\lambda$.}
	\label{fig.CG3D_mesh}%
\end{figure*}

We finally test the potential of our new diffuse interface approach for solving real applications. For this purpose we use a free surface topology based on the real DTM data of the Mont Blanc region\footnote{ The DTM data have been taken from \textcolor{blue}{http://geodati.fmach.it/gfoss\_geodata/libro\_gfoss/.} Our computational domain is centered with respect to the UTM coordinates $(340000.0, 5075000.0)$ }. The horizontal extent of the domain is $28$ km in the $x$ and $y$ directions and ranges from $12$ km below the sea level to $7$ km above it in $z$ direction. We use a heterogeneous material whose parameters are specified in Table \ref{Tab:NT3DCG}. %with $c_p=6000, c_s=3464, \rho=2700$ so that $\lambda=32.4 GPa$ and $\mu=32.4 GPs$ . 
An initial velocity perturbation is placed in $\mathbf{x}_0=(0,0,0)$ for the vertical component of the velocity 
\begin{eqnarray}
	w(\xx,0) = a e^{-r^2/R^2}, 
\label{eq:NTCG3D_1b}
\end{eqnarray}
with $r=\left\| \mathbf{x} \right\|$, $a=-10^{-2}$ and $R=300$m. All other variables for the velocity and the stress tensor are set to zero. The computational domain is covered with a uniform Cartesian grid of $80 \times 80 \times 80$ elements and one refinement level is adopted close to the free surface. In order to represent the complex surface topography within our diffuse interface approach, all that is needed is to compute the shortest distance of a point $\mathbf{x}$ to the free surface defined by the DTM data in order to set the volume fraction function $\alpha$ according to \eqref{eqn.alpha} and \eqref{eq:xi}. The DTM model is given on a Cartesian raster with a spatial resolution of 250 m, which we can then interpolate to any point in our computational domain through bilinear interpolation. 
The smoothing parameter for the diffuse interface zone is set to $I_D=50$ m.  
The simulation with the diffuse interface method is run on the AMR grid with an ADER-DG scheme based on piecewise polinomials of degree $N=3$ in space and time. In Figure \ref{fig.CG3D_mesh} we show a plot of the chosen Cartesian AMR grid with the free surface determined by $\alpha$.

\begin{figure*}
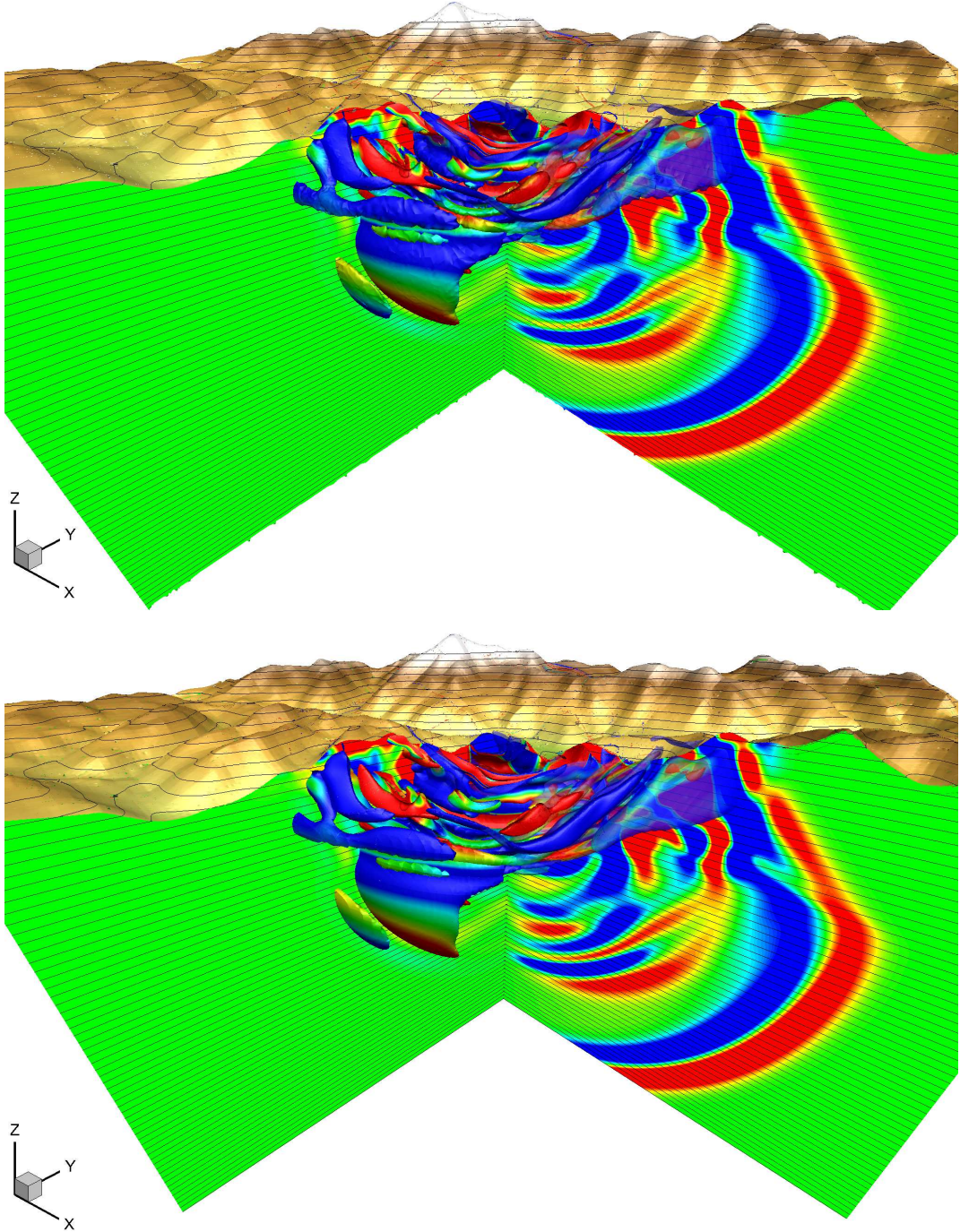

	\centering
\includegraphics[width=0.85\columnwidth]{\flink{}{CG3D_PDESol_t2}} \\
\includegraphics[width=0.85\columnwidth]{\flink{}{CG3D_DILE_t2}}
\caption{Wave propagation in complex 3D geometry. Comparison of the ADER-DG reference solution on unstructured 	
		boundary-fitted grids (left) with the numerical solution obtained with the new diffuse interface method (DIM) on a Cartesian AMR mesh (right) at time $t=2.0$. We show the iso-surfaces $\pm 4\cdot 10^{-5}$ 
		for the velocity components $u$ and $v$ colored by $w$. The slices are colored using the velocity 
		component $u$.}
\label{fig.NTCG1}%
\end{figure*}

\begin{figure*}
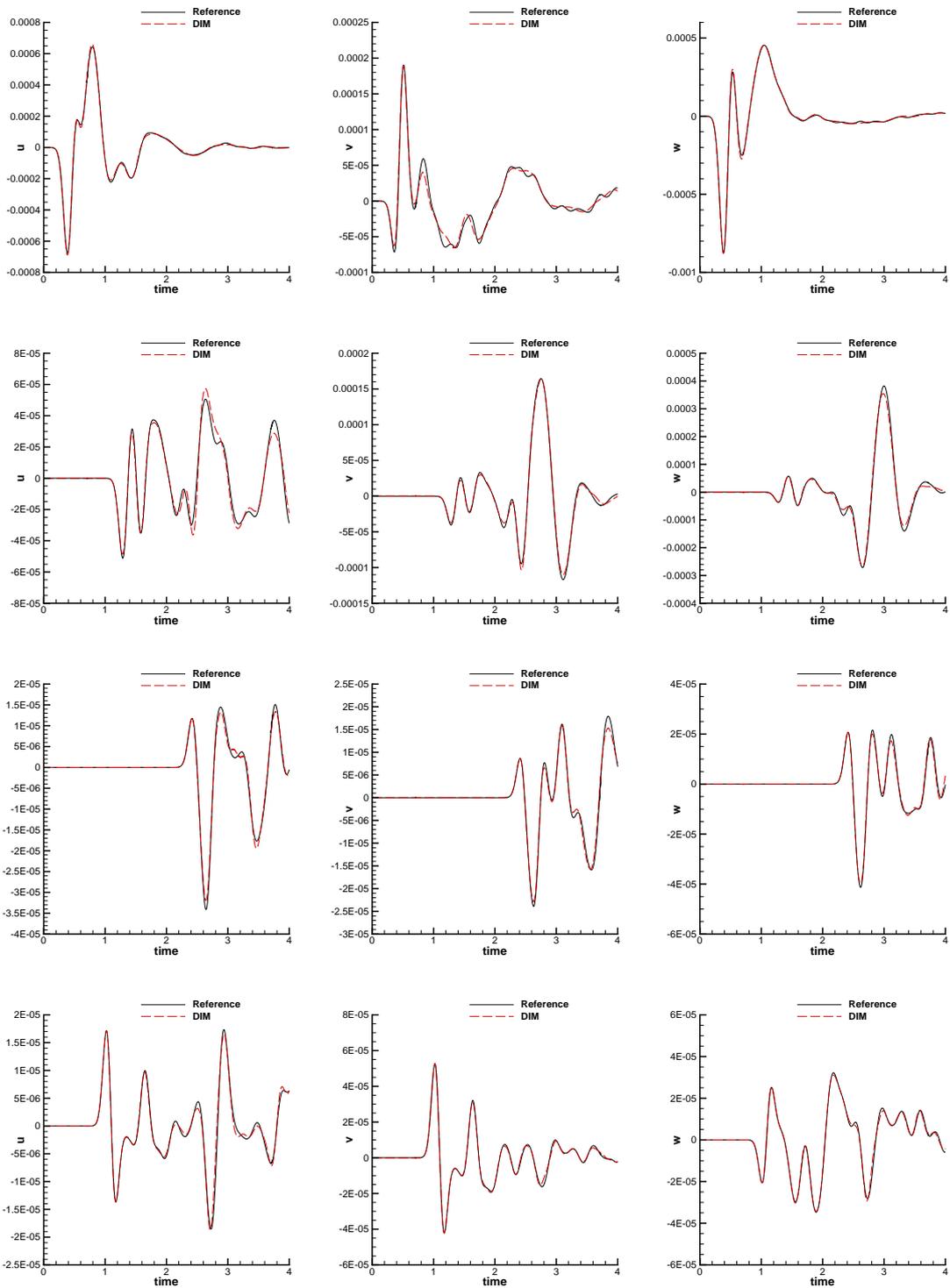
%
	\centering
	\includegraphics[width=0.9\columnwidth]{\flink{}{CG3D_TS1}} 
	\includegraphics[width=0.9\columnwidth]{\flink{}{CG3D_TS2}}
\includegraphics[width=0.9\columnwidth]{\flink{}{CG3D_TS3}} 
\includegraphics[width=0.9\columnwidth]{\flink{}{CG3D_TS4}}
\caption{Wave propagation in complex 3D geometry. Comparison of the time signal of the velocity field obtained with the new diffuse interface approach and the reference solution for the receiver $1$ to $4$ respectively from the top to the bottom row.}
\label{fig.NTCG2}%
\end{figure*}

\begin{figure*}
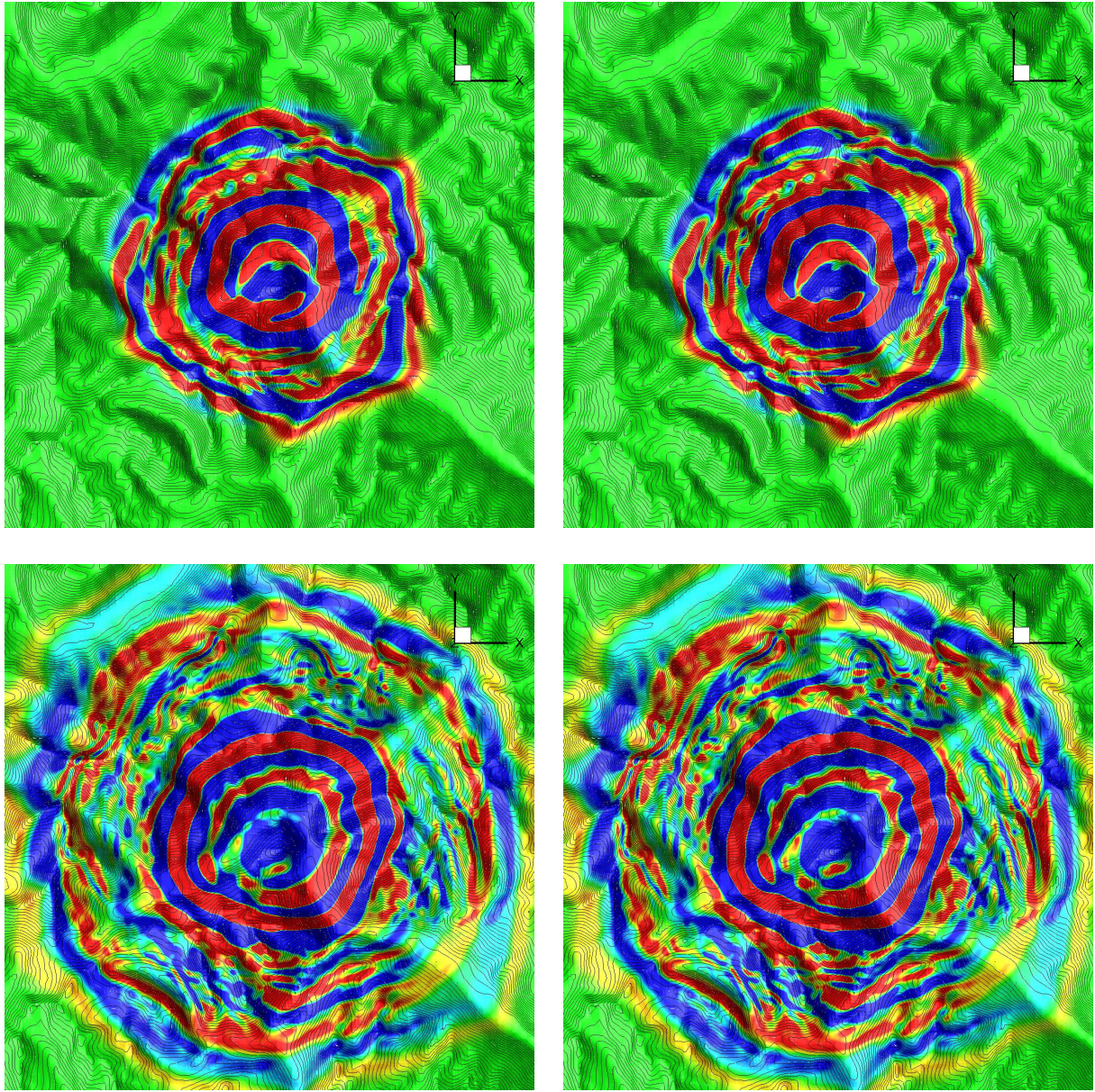
%
	\begin{center}
		\begin{tabular}{lr} 
\includegraphics[width=0.47\columnwidth]{\flink{}{CG3D_DIM_2D_t2}} & 
\includegraphics[width=0.47\columnwidth]{\flink{}{CG3D_PDESOL_2D_t2}} \\ \\ 
\includegraphics[width=0.47\columnwidth]{\flink{}{CG3D_DIM_2D_t3}} & 
\includegraphics[width=0.47\columnwidth]{\flink{}{CG3D_PDESOL_t3}}  
\end{tabular} 
\end{center} 
\caption{Wave propagation in complex 3D geometry. Comparison of the interpolation of the vertical velocity $w$ on the surface computed as the iso-surface $\alpha=0.95$ between the diffuse interface approach (left) and the reference one (right) at time $t=2.0$ and $t=3.0$.}
\label{fig.NTCG3}%
\end{figure*}

The reference solution is computed with an unstructured ADER-DG scheme \cite{gij2} as used in the \texttt{SeisSol} code using $N_e=1,267,717$ boundary-fitted tetrahedral elements and a polynomial approximation degree of $N=3$ in space and time.
A comparison of the numerical solution obtained with the new diffuse interface approach on adaptive Cartesian grids and the results obtained with the unstructured reference code is shown via contour surface maps in Figure \ref{fig.NTCG1} at time $t=2.0$. Overall, we can 
note a very good agreement between the two results.  
We also consider the time signals captured in four receivers, whose positions are reported in Table \ref{tab:NTCG3D_2}. They record data close to the free surface at $1$ km, $5$ km and $10$ km distance from $\mathbf{x}_0$ (receivers $1 \dots 3$) and at $3$ km below the sea level with a distance of $5$ km from $\mathbf{x}_0$ (receiver $4$). The resulting time history of the velocity signals recorded by the four receivers is reported in Figure \ref{fig.NTCG2}. A very good agreement between the new diffuse interface approach and 
the reference scheme can be observed also in this case with complex 3D geometry. Finally, in Figure \ref{fig.NTCG3} we show the interpolation of the velocity component $w$ at the free surface at two different times, where
one can again observe a very good agreement between the numerical results obtained with the new diffuse interface
method and the reference solution obtained on the boundary-fitted unstructured mesh. 

It has to be pointed out that the setup of this test problem with the new diffuse interface approach is  \textit{completely automatic}, without requiring \textit{any} manual user interaction. The entire setup 
process of the computational model starts with reading the DTM data from a file according  
to well-established standard GIS file formats, continues by automatically setting the color function $\alpha$  according to \eqref{eqn.alpha} and \eqref{eq:xi} with appropriate bilinear interpolation of the DTM data 
to the nodal degrees of freedom of the ADER-DG scheme and to the subcell FV averages and closes with the 
automatic setup of the adaptive Cartesian AMR grid based on the gradient of $\alpha$ up to the desired level of spatial resolution. We would like to emphasize again that for the diffuse 
interface approach the time step size does \textit{not} depend on the distribution of $\alpha$. 
In contrast to this fully automated chain in \texttt{ExaHyPE}, the setup of the same test case in \texttt{SeisSol}  
still requires the generation of a boundary-aligned unstructured tetrahedral mesh with an external grid 
generation tool that needs some manual interactions with the end user. For very complex
surface topography, even more user interaction is required to obtain a high quality grid, which is essential 
due to the CFL restriction on the time step. 

\begin{table*}%
	\begin{center}
		\begin{tabular}{ccccccc}
			& Position 	& $c_p (ms^{-1})$ & $c_s (ms^{-1})$ & $\rho (kg m^{-3})$ & $\lambda (GPa)$ & $\mu (GPa)$	 \\
			\hline
			Medium 1	& $z>-1000$m		  &	4000	&  2000 & 2600 & 20.8 & 10.4 \\
			Medium 2	& $z\leq-1000$m		&	6000	&  3464 & 2700 & 32.4 & 32.4 \\
			\hline
		\end{tabular}
	\end{center}
	\caption{Material parameters for the wave propagation test in a complex 3D geometry.}
	\label{Tab:NT3DCG}
\end{table*}

\begin{table}%
\begin{center}
\begin{tabular}{cccc}
	Receiver & $x$ 	& $y$	 & $z$ \\
	\hline
1&1000.000000 & 0.000000 &1397.723250  \\ % 9, 1km from Epicenter
2&3535.533906 &3535.533906 &1883.989778 \\ %16 5km from Epicenter
3&8660.254038 &5000.000000 &2173.363299 \\ %24 10km from Epicenter
4&1545.084972 &4755.282581 &-3000.000000 \\ % 41, deep, 5km from epicenter
	\hline
\end{tabular}
\end{center}
\caption{Receiver positions for the wave propagation test in complex 3D geometry.} 
\label{tab:NTCG3D_2}
\end{table}

\section{Conclusions}
\label{sec:conclusions}

In this paper a novel diffuse interface method (DIM) for the simulation of seismic wave propagation in linear  isotropic  material with complex free surface topography has been proposed. The governing PDE system can be derived from the Baer-Nunziato (BN) model of compressible multi-phase  flows \cite{BaerNunziato1986,SaurelAbgrall,SaurelAbgrall2,AbgrallSaurel} following similar ideas as those employed in \cite{DIM2D,DIM3D}.      
In alternative, our governing equations can also be derived by combining the equations of nonlinear hyperelasticity in Eulerian coordinates of Godunov and Romenski \cite{GodunovRomenski72,Godunov:2003a,Godunov:1995a} with the compressible multi-phase flow model of Romenski et al. \cite{Rom1998,Rom2010}. 

In both cases, the velocity of the medium is supposed to be very small, so that the nonlinear convective terms can be neglected, and a linear material behavior according to Hooke's law is assumed. We have proven that the solution  of the Riemann problem with arbitrary data and a 
jump of the volume fraction function $\alpha$ from one to zero yields a Godunov state at the  interface in which the normal components of the stress tensor vanish, which  
is exactly the required free surface boundary condition $\boldsymbol{\sigma} \cdot \mathbf{n} = 0$. The free surface boundary condition on a physical domain $\Omega_p$ of arbitrary shape can therefore be properly imposed by defining a spatially  variable scalar function $\alpha$ in the computational domain $\Omega$, which has to be large enough to contain $\Omega_p$, i.e.~$\Omega_p \subset \Omega$, simply by setting $\alpha=1$ for $\mathbf{x} \in \Omega_p$ and $\alpha=0$ for $\mathbf{x} \notin \Omega_p$, 
without having to fit the boundary of the computational domain $\partial \Omega$ to the real boundary $\partial \Omega_p$ of the physical
domain to be discretized. In practical simulations, the interface layer which contains the transition from $\alpha=1$ to $\alpha=0$ is slightly smoothed by a characteristic width $I_D$, which is the reason why we call our approach a diffuse interface method.  
We have carried out a systematic study in which we show that for vanishing interface thickness $I_D \to 0$ the correct wave reflection is obtained.  

The governing equations derived in the first part of the paper have been solved on  adaptive Cartesian meshes (AMR) via high order accurate ADER-DG schemes  combined with a sub-cell finite volume limiter \cite{Dumbser2014,Zanotti2015}. The use
of the subcell finite volume limiter is necessary in regions with strong gradients of $\alpha$
in order to avoid spurious oscillations and unphysical solutions that would be obtained with 
a pure unlimited high order DG scheme.    
The practical implementation of the model  has been carried out in the \texttt{ExaHyPE} code developed within the European H2020 research project \textit{An Exascale Hyperbolic PDE Engine}, see \textcolor{blue}{http://exahype.eu/}. We have presented a large set of two- and three-dimensional wave propagation problems where we have compared the results obtained with the new diffuse interface approach with classical computational methods based on boundary-fitted unstructured meshes. In all cases under investigation, and even in  the presence of complex surface topography, the new diffuse interface model performs very well as shown throughout this paper. 

We stress again that the \textit{key novelty} introduced here consists in the representation of the  geometrically complex surface topography merely via the scalar solid volume fraction function $\alpha$, instead of making use of complex structured or unstructured boundary-fitted meshes. In order to improve the spatial resolution of certain geometric features of the physical domain $\Omega_p$, we simply use adaptive mesh refinement (AMR) on locally Cartesian grids. This allows a \textit{fully automated workflow} in the setup of the computational model, without requiring any external mesh generation tools or any manual interaction with the user. We underline again that the time step restriction in our new approach is completely independent of the complexity of the geometry of the domain $\Omega_p$ to be discretized, since $\alpha$ has no influence on the 
eigenvalues of the governing PDE system. The admissible local time step size according to the CFL condition is therefore only given by the local mesh size $h$, the pressure wave propagation speed 
$c_p$ and the polynomial approximation degree $N$. 

Current work in progress is the implementation of new strategies for highly efficient   
small matrix-matrix multiplications in ADER-DG schemes on adaptive Cartesian grids 
(exploiting also the fact that we use a nodal tensor-product basis) in order to improve computational 
performance of the code, similar to the hardware optimizations already successfully applied 
in the context of the unstructured ADER-DG schemes used in \texttt{SeisSol} \cite{SeisSol1,SeisSol2}.  

Future research will concern the extension of our new diffuse interface approach to the full equations
of nonlinear hyperelasticity including plastic deformations and dynamic rupture processes based on the  
Godunov-Peshkov-Romenski model presented and discussed in  \cite{GodunovRomenski72,Godunov:2003a,Godunov:1995a,PeshRom2014,DumbserPeshkov2016,DumbserPeshkov2017}.       

%-------------------------------------------------------------
\section*{Acknowledgements}
%-------------------------------------------------------------

This research was funded by the European Union's Horizon 2020 
Research and Innovation Programme under the project \textit{ExaHyPE},
grant no. 671698 (call FETHPC-1-2014). The 3D simulations were performed 
on the HazelHen supercomputer at the HLRS in Stuttgart, Germany and 
on the SuperMUC supercomputer at the LRZ in Garching, Germany.

\bibliographystyle{elsarticle-num}
\bibliography{AMRDILE}

\begin{thebibliography}{100}
\expandafter\ifx\csname url\endcsname\relax
  \def\url#1{\texttt{#1}}\fi
\expandafter\ifx\csname urlprefix\endcsname\relax\def\urlprefix{URL }\fi
\expandafter\ifx\csname href\endcsname\relax
  \def\href#1#2{#2} \def\path#1{#1}\fi

\bibitem{Madariaga1976}
R.~Madariaga, Dynamics of an expanding circular fault, Bulletin of the
  Seismological Society of America 66 (1976) 639--666.

\bibitem{Virieux1984}
J.~Virieux, Sh-wave propagation in heterogeneous media: Velocity--stress
  finite--difference method, Geophysics 49 (1984) 1933--1942.

\bibitem{Virieux1986}
J.~Virieux, P-sv wave propagation in heterogeneous media: Velocity--stress
  finite--difference method, Geophysics 51 (1986) 889--901.

\bibitem{Levander1988}
A.~Levander, Fourth-order finite difference p-sv seismograms, Geophysics 53
  (1988) 1425--1436.

\bibitem{Mora1989}
P.~Mora, Modeling anisotropic seismic waves in 3-d, SEG society of exploration
  Geophysicists (1989) 1039--1043.

\bibitem{moczo2002}
P.~Moczo, J.~Kristek, V.~Vavrycuk, R.~Archuleta, L.~Halada, {3D heterogeneous
  staggered-grid finite-difference modeling of seismic motion with volume
  harmonic and arithmetic averaging of elastic moduli and densities}, Bulletin
  of the Seismological Society of America 92 (2002) 3042--3066.

\bibitem{igel1995}
H.~Igel, P.~Mora, B.~Riollet, Anisotropic wave propagation through
  finite-difference grids, Geophysics 60 (1995) 1203--1216.

\bibitem{tessmer1995}
E.~Tessmer, 3-d seismic modelling of general material anisotropy in the
  presence of the free surface by a chebyshev spectral method, Geophysical
  Journal International 121 (1995) 557--575.

\bibitem{magnier1994}
S.~Magnier, P.~Mora, Finite differences on minimal grids, Geophysics 59 (1994)
  1435--1443.

\bibitem{kaser2001a}
M.~K\"aser, H.~Igel, Numerical simulation of 2d wave propagation on
  unstructured grids using explicit differential operators, Geophysical
  Prospecting 49 (2001) 607--619.

\bibitem{kaser2001b}
M.~K\"aser, H.~Igel, A comparative study of explicit differential operators on
  arbitrary grids, Journal of Computational Acoustics 9 (2001) 1111--1125.

\bibitem{Wang2002a}
Z.~Wang, {Spectral (finite) volume method for conservation laws on unstructured
  grids: basic formulation}, Journal of Computational Physics 178 (2002)
  210--251.

\bibitem{Wang2002b}
Z.~Wang, Y.~Liu, {Spectral (finite) volume method for conservation laws on
  unstructured grids II: extension to two-dimensional scalar equation}, Journal
  of Computational Physics 179 (2002) 665--697.

\bibitem{Wang2004a}
Z.~Wang, Y.~Liu, {Spectral (finite) volume method for conservation laws on
  unstructured grids III: one-dimensional systems and partition optimization},
  Journal of Scientific Computing 20 (2004) 137--157.

\bibitem{Wang2004b}
Z.~Wang, Y.~Liu, {Spectral (finite) volume method for conservation laws on
  unstructured grids IV: extension to two-dimensional systems}, Journal of
  Computational Physics 194 (2004) 716--741.

\bibitem{Tadi2004}
M.~Tadi, {Finite Volume Method for 2D Elastic Wave Propagation}, Bulletin of
  the Seismological Society of America 94 (2004) 1500--1509.

\bibitem{Dormt1995}
E.~Dormy, T.~A, {Numerical simulation of elastic wave propagation using a
  finite volume method}, Journal of Geophysical research 100 (1995) 2123--2133.

\bibitem{dumbserkaeser06c}
M.~Dumbser, M.~K\"aser, Arbitrary high order finite volume schemes for seismic
  wave propagation on unstructured meshes in 2d and 3d, Geophysical Journal
  International 171 (2007) 665--694.

\bibitem{Moczo2006}
M.~Kristekov\'a, J.~Kristek, P.~Moczo, S.~Day, {Misfit criteria for
  quantitative comparison of seismograms}, Bulletin of the Seismological
  Society of America 96 (2006) 1836--1850.

\bibitem{Moczo2009}
M.~Kristekov\'a, J.~Kristek, P.~Moczo, {Time-frequency misfit and
  goodness-of-fit criteria for quantitative comparison of time signals},
  Geophysical Journal International 178 (2009) 813--825.

\bibitem{Hermann2008}
M.~K\"aser, V.~Hermann, J.~de~la Puente, {Quantitative accuracy analysis of the
  discontinuous Galerkin method for seismic wave propagation}, Geophysical
  Journal International 173 (2008) 990--999.

\bibitem{Moczo2010a}
P.~Moczo, J.~Kristek, M.~Galis, P.~Pazak, {On accuracy of the finite-difference
  and finite-element schemes with respect to P-wave to S-wave speed ratio},
  Geophysical Journal International 182 (2010) 493--510.

\bibitem{Patera1984}
A.~T. Patera, A spectral-element method for fluid dynamics: laminar flow in a
  channel expansion, Journal of Computational Physics 144 (1984) 45--58.

\bibitem{Priolo1994}
E.~Priolo, J.~Carcione, G.~Seriani, Numerical simulation of interface waves by
  high-order spectral modeling techniques, Journal of Computational Physics 144
  (1984) 45--58.

\bibitem{Komatitsch1998}
D.~Komatitsch, J.~Vilotte, The spectral-element method: an efficient tool to
  simulate the seismic response of 2d and 3d geological structures, Bulletin of
  the Seismological Society of America 88 (1998) 368--392.

\bibitem{Seriani1998}
G.~Seriani, 3-d large-scale wave propagation modeling by a spectral-element
  method on a cray t3e multiprocessor, Computer Methods in Applied Mechanics
  and Engineering 164 (1998) 235--247.

\bibitem{Komatitsch1999}
D.~Komatitsch, J.~Tromp, Introduction to the spectral-element method for 3-d
  seismic wave propagation, Geophysical Journal International 139 (1999)
  806--822.

\bibitem{Komatitsch2002}
D.~Komatitsch, J.~Tromp, Spectral-element simulations of global seismic wave
  propagation—i. validation, Geophysical Journal International 149 (2002)
  390--412.

\bibitem{Tessmer1994}
E.~Tessmer, D.~Kosloff, 3-d elastic modelling with surface topography by a
  chebyshev spectral method, Geophysics 59 (1994) 464--473.

\bibitem{Igel1999}
H.~Igel, Wave propagation in three-dimensional spherical sections by the
  chebyshev spectral method, Geophysical Journal International 136 (1999)
  559--566.

\bibitem{gugubrian_waves}
G.~Scovazzi, B.~Carnes, Weak boundary conditions for wave propagation problems
  in confined domains: Formulation and implementation using a variational
  multiscale method, Computer Methods in Applied Mechanics and Engineering
  221-222 (2012) 117--131.

\bibitem{song2015nitsche}
T.~Song, G.~Scovazzi, {A Nitsche method for wave propagation problems in time
  domain}, Computer Methods in Applied Mechanics and Engineering 293 (2015)
  481--521.

\bibitem{scovazzi2017velocity}
G.~Scovazzi, T.~Song, X.~Zeng, A velocity/stress mixed stabilized nodal finite
  element for elastodynamics: Analysis and computations with strongly and
  weakly enforced boundary conditions, Computer Methods in Applied Mechanics
  and Engineering 325 (2017) 532--576.

\bibitem{gij1}
M.~K\"aser, M.~Dumbser, {An arbitrary high-order discontinuous Galerkin method
  for elastic waves on unstructured meshes – I. The two-dimensional isotropic
  case with external source terms}, Geophysical Journal International 166
  (2006) 855--877.

\bibitem{gij2}
M.~Dumbser, M.~K\"aser, {An arbitrary high-order discontinuous Galerkin method
  for elastic waves on unstructured meshes – II. The three-dimensional
  isotropic case}, Geophysical Journal International 167 (2006) 319--336.

\bibitem{gij5}
M.~Dumbser, M.~K\"aser, E.~F. Toro, {An arbitrary high--order Discontinuous
  Galerkin method for elastic waves on unstructured meshes -- V. Local time
  stepping and p-adaptivity}, Geophysical Journal International 171 (2007)
  695--717.

\bibitem{GroteDG}
M.~Grote, A.~Schneebeli, D.~Sch\"otzau, {Discontinuous Galerkin finite element
  method for the wave equation}, SIAM Journal on Numerical Analysis 44 (2006)
  2408--2431.

\bibitem{Antonietti1}
P.~Antonietti, I.~Mazzieri, A.~Quarteroni, F.~Rapetti, {Non-conforming high
  order approximations of the elastodynamics equation}, Computer Methods in
  Applied Mechanics and Engineering 209--212 (2012) 212--238.

\bibitem{Antonietti2}
P.~Antonietti, C.~Marcati, I.~Mazzieri, A.~Quarteroni, {High order
  discontinuous Galerkin methods on simplicial elements for the elastodynamics
  equation}, Numerical Algorithms 71 (2016) 181--206.

\bibitem{spacetimedg1}
J.~J.~W. van~der Vegt, H.~van~der Ven, Space–-time discontinuous {Galerkin}
  finite element method with dynamic grid motion for inviscid compressible
  flows {I}. general formulation, Journal of Computational Physics 182 (2002)
  546–--585.

\bibitem{spacetimedg2}
H.~van~der Ven, J.~J.~W. van~der Vegt, Space-–time discontinuous {Galerkin}
  finite element method with dynamic grid motion for inviscid compressible
  flows {II}. efficient flux quadrature, Comput. Methods Appl. Mech. Engrg. 191
  (2002) 4747–--4780.

\bibitem{KlaijVanDerVegt}
C.~Klaij, J.~J. W.~V. der Vegt, H.~V. der Ven, {Space-time discontinuous
  Galerkin method for the compressible Navier-Stokes equations}, Journal of
  Computational Physics 217 (2006) 589--611.

\bibitem{Rhebergen2012}
S.~Rhebergen, B.~Cockburn, {A space–time hybridizable discontinuous Galerkin
  method for incompressible flows on deforming domains}, Journal of
  Computational Physics 231 (2012) 4185--4204.

\bibitem{Rhebergen2013}
S.~Rhebergen, B.~Cockburn, J.~J. van~der Vegt, {A space-time discontinuous
  Galerkin method for the incompressible Navier-Stokes equations}, Journal of
  Computational Physics 233 (2013) 339--358.

\bibitem{Balazsova1}
M.~Balazsova, M.~Feistauer, {On the stability of the ALE space-time
  discontinuous Galerkin method for nonlinear convection-diffusion problems in
  time-dependent domains}, Applications of Mathematics 60 (2015) 501--526.

\bibitem{Balazsova2}
M.~Balazsova, M.~Feistauer, M.~Hadrava, A.~Kosik, {On the stability of the
  space-time discontinuous Galerkin method for the numerical solution of
  nonstationary nonlinear convection-diffusion problems}, Journal of Numerical
  Mathematics 23 (2015) 211--233.

\bibitem{Antonietti3}
P.~Antonietti, I.~Mazzieri, A.~Quarteroni, F.~Rapetti, {High order space-time
  discretization for elastic wave propagation problems}, in: M.~Azaiez, H.~E.
  Fekihand, J.~Hestaven (Eds.), Proceedings of ICOSAHOM 2012, LNCSE, Vol.~95,
  Springer Verlag, 2014, pp. 87--97.

\bibitem{Antonietti4}
P.~Antonietti, N.~D. Santo, I.~Mazzieri, A.~Quarteroni, {A high-order
  discontinuous Galerkin approximation to ordinary differential equations with
  applications to elastodynamics}, IMA Journal of Numerical Analysis.

\bibitem{SIStag2017}
M.~Tavelli, M.~Dumbser, {Arbitrary high order accurate space-time discontinuous
  Galerkin finite element schemes on staggered unstructured meshes for linear
  elasticity}, Journal of Computational Physics {366} (2018) {386--414}.

\bibitem{DumbserCasulli}
M.~Dumbser, V.~Casulli, A staggered semi-implicit spectral discontinuous
  {Galerkin} scheme for the shallow water equations, Applied Mathematics and
  Computation 219~(15) (2013) 8057--8077.

\bibitem{2DSIUSW}
M.~Tavelli, M.~Dumbser, {A high order semi-implicit discontinuous Galerkin
  method for the two dimensional shallow water equations on staggered
  unstructured meshes}, Applied Mathematics and Computation 234 (2014)
  623--644.

\bibitem{2STINS}
M.~Tavelli, M.~Dumbser, {A staggered arbitrary high order semi-implicit
  discontinuous Galerkin method for the two dimensional incompressible
  Navier-Stokes equations}, Computers and Fluids 119 (2015) 235--249.

\bibitem{3DSIINS}
M.~Tavelli, M.~Dumbser, {A staggered, space-time discontinuous Galerkin method
  for the three-dimensional incompressible Navier-Stokes equations on
  unstructured tetrahedral meshes}, Journal of Computational Physics 319 (2016)
  294--323.

\bibitem{Fambri2016}
F.~Fambri, M.~Dumbser, {Spectral semi-implicit and space-time discontinuous
  Galerkin methods for the incompressible Navier-Stokes equations on staggered
  Cartesian grids}, Applied Numerical Mathematics 110 (2016) 41--74.

\bibitem{3DSICNS}
M.~Tavelli, M.~Dumbser, {A pressure--based semi--implicit space–time
  discontinuous Galerkin method on staggered unstructured meshes for the
  solution of the compressible Navier--Stokes equations at all Mach numbers},
  Journal of Computational Physics 341 (2017) 341--376.

\bibitem{AMRDGSI}
F.~Fambri, M.~Dumbser, Semi-implicit discontinuous {G}alerkin methods for the
  incompressible {N}avier-{S}tokes equations on adaptive staggered {C}artesian
  grids, Computer Methods in Applied Mechanics and Engineering 324 (2017)
  170--203.

\bibitem{Bern1992}
M.~Bern, D.~Eppstein, Mesh generation and optimal triangulation, Computing in
  Euclidean Geometry 1 (1992) 23--90.

\bibitem{Joe1995}
B.~Joe, Construction of three-dimensional improved-quality triangulations using
  local transformations, SIAM Journal of Scientific COmputing 16 (1995)
  1292--1307.

\bibitem{Fleischmann1999}
P.~Fleischmann, W.~Pyka, S.~Selberherr, Mesh generation for application in
  technology cad, IEICE Transactions on Electronics E82--C (1999) 937--947.

\bibitem{Cheng2000}
S.~W. Cheng, T.~K. Dey, H.~Edelsbrunner, , M.~A. Facello, S.~H. Teng, Sliver
  exudation, Journal of the ACM 47 (2000) 883--904.

\bibitem{Edelsbrunner2002}
H.~Edelsbrunner, D.~Guoy, An experimental study of sliver exudation,
  Engineering with Computers 18 (2002) 229--240.

\bibitem{Taub2009}
A.~Taube, M.~Dumbser, C.~Munz, R.~Schneider, {A High Order Discontinuous
  Galerkin Method with Local Time Stepping for the Maxwell Equations},
  International Journal Of Numerical Modelling: Electronic Networks, Devices
  And Fields 22 (2009) 77--103.

\bibitem{GroteLTS1}
M.~Grote, T.~Mitkova, {High-order explicit local time-stepping methods for
  damped wave equations}, Journal of Computational and Applied Mathematics 239
  (2013) 270--289.

\bibitem{GroteLTS2}
M.~Grote, T.~Mitkova, {Explicit local time-stepping methods for Maxwell's
  equations}, Journal of Computational and Applied Mathematics 234 (2010)
  3283--3302.

\bibitem{BaerNunziato1986}
M.~R. Baer, J.~W. Nunziato, A two-phase mixture theory for the
  deflagration-to-detonation transition {(DDT)} in reactive granular materials,
  J. Multiphase Flow 12 (1986) 861--889.

\bibitem{SaurelAbgrall}
R.~Saurel, R.~Abgrall, {A Multiphase Godunov Method for Compressible Multifluid
  and Multiphase Flows}, Journal of Computational Physics 150 (1999) 425--467.

\bibitem{SaurelAbgrall2}
R.~Saurel, R.~Abgrall, {A Simple Method for Compressible Multifluid Flows},
  SIAM Journal on Scientific Computing 21 (1999) 1115--1145.

\bibitem{Abgral2003}
R.~Abgrall, B.~Nkonga, R.~Saurel, Efficient numerical approximation of
  compressible multi-material flow for unstructured meshes, Computers and
  Fluids 32 (2003) 571--605.

\bibitem{AbgrallSaurel}
R.~Abgrall, R.~Saurel, Discrete equations for physical and numerical
  compressible multiphase mixtures, Journal of Computational Physics 186 (2003)
  361--396.

\bibitem{GodunovRomenski72}
S.~K. Godunov, E.~I. Romenski, Nonstationary equations of the nonlinear theory
  of elasticity in {Euler} coordinates., Journal of Applied Mechanics and
  Technical Physics 13 (1972) 868--885.

\bibitem{Rom1998}
E.~Romenski, Hyperbolic systems of thermodynamically compatible conservation
  laws in continuum mechanics, Mathematical and computer modelling 28 (1998)
  115--130.

\bibitem{PeshRom2014}
I.~Peshkov, E.~Romenski, A hyperbolic model for viscous {{N}ewtonian} flows,
  Continuum Mechanics and Thermodynamics 28 (2016) 85--104.

\bibitem{DumbserPeshkov2016}
M.~Dumbser, I.~Peshkov, E.~Romenski, O.~Zanotti, {High order ADER schemes for a
  unified first order hyperbolic formulation of continuum mechanics: Viscous
  heat-conducting fluids and elastic solids}, Journal of Computational Physics
  314 (2016) 824--862.

\bibitem{DumbserPeshkov2017}
M.~Dumbser, I.~Peshkov, E.~Romenski, O.~Zanotti, {High order ADER schemes for a
  unified first order hyperbolic formulation of Newtonian continuum mechanics
  coupled with electro-dynamics}, Journal of Computational Physics 348 (2017)
  298--342.

\bibitem{Romenskii2017}
E.~I. Romenskii, E.~B. Lys', V.~A. Cheverda, M.~I. Epov, {Dynamics of
  deformation of an elastic medium with initial stresses}, Journal of Applied
  Mechanics and Technical Physics 58 (2017) 914--923.

\bibitem{Ndanou2015}
S.~Ndanou, N.~Favrie, S.~Gavrilyuk, Multi-solid and multi-fluid diffuse
  interface model: Applications to dynamic fracture and fragmentation, Journal
  of Computational Physics 295 (2015) 523--555.

\bibitem{FavrieGavrilyuk2012}
N.~Favrie, S.~L. Gavrilyuk, Diffuse interface model for compressible
  fluid-compressible elastic-plastic solid interaction, Journal of
  Computational Physics 231 (2012) 2695--2723.

\bibitem{FavrieGavrilyuk2009}
N.~Favrie, S.~L. Gavrilyuk, Solid-fluid diffuse interface model in cases of
  extreme deformations, Journal of Computational Physics 228 (2009) 6037--6077.

\bibitem{Gavrilyuk2008}
S.~L. Gavrilyuk, N.~Favrie, R.~Saurel, Modelling wave dynamics of compressible
  elastic materials, Journal of Computational Physics 227 (2008) 2941--2969.

\bibitem{Saurel2009}
R.~Saurel, F.~Petitpas, R.~Berry, Simple and efficient relaxation method for
  interfaces separating compressible fluids cavitating flows and shock in
  multiphase mixtures, Journal of Computational Physics 228 (2009) 1678--1712.

\bibitem{Kapila2001}
A.~Kapila, R.~Menikoff, J.~Bdzil, S.~Son, D.~Stewart, Two-phase modeling of ddt
  in granular materials: reduced equations, Physics of Fluids 13 (2001)
  3002--3024.

\bibitem{Zanotti2015}
O.~Zanotti, F.~Fambri, M.~Dumbser, A.~Hidalgo, {Space-time adaptive ADER
  discontinuous Galerkin finite element schemes with a posteriori sub-cell
  finite volume limiting}, Computers and Fluids 118 (2015) 204--224.

\bibitem{Dumbser2014}
M.~Dumbser, O.~Zanotti, R.~Loub\`ere, S.~Diot, {A posteriori subcell limiting
  of the discontinuous Galerkin finite element method for hyperbolic
  conservation laws}, Journal of Computational Physics 278 (2014) 47--75.

\bibitem{Clain2011}
S.~Clain, S.~Diot, R.~Loub\`ere, A high-order finite volume method for systems
  of conservation lawsmulti-dimensional optimal order detection (mood), Journal
  of Computational Physics 230 (2011) 4028--4050.

\bibitem{Diot2012}
S.~Diot, S.~Clain, R.~Loub\`ere, Improved detection criteria for the
  multi-dimensional optimal order detection (mood) on unstructured meshes with
  very high-order polynomials, Computers and Fluids 64 (2012) 43--63.

\bibitem{HLLEM}
M.~Dumbser, D.~Balsara, {A new efficient formulation of the HLLEM Riemann
  solver for general conservative and non-conservative hyperbolic systems},
  Journal of Computational Physics 304 (2016) 275--319.

\bibitem{NCP_HLLEM}
M.~Dumbser, D.~S. Balsara, A new efficient formulation of the \{HLLEM\} riemann
  solver for general conservative and non-conservative hyperbolic systems,
  Journal of Computational Physics 304 (2016) 275--319.

\bibitem{BedfordDrumheller}
A.~Bedford, D.~Drumheller, Elastic Wave Propagation, Wiley, Chichester, UK,
  1994.

\bibitem{DIM2D}
M.~Dumbser, A simple two-phase method for the simulation of complex free
  surface flows, Computer Methods in Applied Mechanics and Engineering 200
  (2011) 1204--1219.

\bibitem{DIM3D}
M.~Dumbser, {A Diffuse Interface Method for Complex Three-Dimensional Free
  Surface Flows}, Computer Methods in Applied Mechanics and Engineering 257
  (2013) 47--64.

\bibitem{Godunov:2003a}
S.~K. Godunov, E.~I. Romenski, {Elements of Continuum Mechanics and
  Conservation Laws}, Kluwer Academic/ Plenum Publishers, 2003.

\bibitem{Godunov:1995a}
S.~K. Godunov, E.~I. Romenski, {Thermodynamics, conservation laws, and
  symmetric forms of differential equations in mechanics of continuous media},
  in: {Computational Fluid Dynamics Review 95}, {John Wiley, NY}, 1995, pp.
  19--31.

\bibitem{Rom2010}
E.~Romenski, D.~Drikakis, E.~Toro, Conservative models and numerical methods
  for compressible two-phase flow, Journal of Scientific Computing 42 (2010)
  68--95.

\bibitem{Castro2006}
M.~J. Castro, J.~M. Gallardo, C.~Par\'es, High-order finite volume schemes
  based on reconstruction of states for solving hyperbolic systems with
  nonconservative products. applications to shallow-water systems, Mathematics
  of Computations 75 (2006) 1103--1134.

\bibitem{Pares2006}
C.~Par\'es, Numerical methods for nonconservative hyperbolic systems: a
  theoretical framework, SIAM Journal on Numerical Analysis 44 (2006) 300--321.

\bibitem{OsherNC}
M.~Dumbser, E.~F. Toro, A simple extension of the {Osher} {Riemann} solver to
  non-conservative hyperbolic systems, Journal of Scientific Computing 48
  (2011) 70--88.

\bibitem{DLMtheory}
G.~D. Maso, P.~G. LeFloch, F.~Murat, Definition and weak stability of
  nonconservative products, J. Math. Pures Appl. 74 (1995) 483--548.

\bibitem{Dumbser2008}
M.~Dumbser, D.~S. Balsara, E.~F. Toro, C.~D. Munz, A unified framework for the
  construction of one-step finite-volume and discontinuous {Galerkin} schemes,
  Journal of Computational Physics 227 (2008) 8209--–8253.

\bibitem{AMR3DCL}
M.~Dumbser, O.~Zanotti, A.~Hidalgo, D.~Balsara, {ADER-{WENO} Finite Volume
  Schemes with Space-Time Adaptive Mesh Refinement}, Journal of Computational
  Physics 248 (2013) 257--286.

\bibitem{Peano1}
H.~Bungartz, M.~Mehl, T.~Neckel, T.~Weinzierl, {The PDE framework Peano applied
  to fluid dynamics: An efficient implementation of a parallel multiscale fluid
  dynamics solver on octree-like adaptive Cartesian grids}, Computational
  Mechanics 46 (2010) 103--114.

\bibitem{Peano2}
T.~Weinzierl, M.~Mehl, {Peano-A traversal and storage scheme for octree-like
  adaptive Cartesian multiscale grids}, SIAM Journal on Scientific Computing 33
  (2011) 2732--2760.

\bibitem{Khokhlov1998}
A.~Khokhlov, Fully threaded tree algorithms for adaptive refinement fluid
  dynamics simulations, Journal of Computational Physics 143~(2) (1998) 519 --
  543.

\bibitem{Zanotti2015c}
O.~Zanotti, F.~Fambri, M.~Dumbser, A.~Hidalgo, Space-time adaptive {ADER}
  discontinuous {{G}alerkin} finite element schemes with a posteriori sub-cell
  finite volume limiting, Computers and Fluids 118 (2015) 204 -- 224.

\bibitem{Zanotti2015d}
O.~Zanotti, F.~Fambri, M.~Dumbser, Solving the relativistic
  magnetohydrodynamics equations with {ADER} discontinuous {G}alerkin methods,
  a posteriori subcell limiting and adaptive mesh refinement, Mon. Not. R.
  Astron. Soc. 452 ({2015}) 3010--3029.

\bibitem{ADERDGVisc}
F.~Fambri, M.~Dumbser, O.~Zanotti, Space-time adaptive {ADER}-{DG} schemes for
  dissipative flows: {C}ompressible {N}avier-{S}tokes and resistive {MHD}
  equations, Computer Physics Communications 220 (2017) 297--318.

\bibitem{gji1}
M.~K\"aser, M.~Dumbser, {An arbitrary high-order discontinuous Galerkin method
  for elastic waves on unstructured meshes - I. The two-dimensional isotropic
  case with external source terms.}, Geophysical Journal International 166
  (2006) 855--877.

\bibitem{AMR3DNC}
M.~Dumbser, A.~Hidalgo, O.~Zanotti, {High Order Space-Time Adaptive ADER-WENO
  Finite Volume Schemes for Non-Conservative Hyperbolic Systems}, Computer
  Methods in Applied Mechanics and Engineering 268 (2014) 359--387.

\bibitem{SeisSol1}
A.~Breuer, A.~Heinecke, M.~Bader, C.~Pelties, {Accelerating SeisSol by
  generating vectorized code for sparse matrix operators}, Advances in Parallel
  Computing 25 (2014) 347--356.

\bibitem{SeisSol2}
A.~Breuer, A.~Heinecke, S.~Rettenberger, M.~Bader, A.~Gabriel, C.~Pelties,
  {Sustained petascale performance of seismic simulations with SeisSol on
  SuperMUC}, Lecture Notes in Computer Science (LNCS) 8488 (2014) 1--18.

\end{thebibliography}
\end{document}